\documentclass[11pt]{article}

\newtheorem{theorem}{Theorem}

\newtheorem{claim}{Claim}
\newtheorem{proposition}{Proposition}
\newtheorem{conjecture}{Conjecture}
\newtheorem{lemma}{Lemma}

\usepackage{amsmath,amssymb,color,graphicx,epsf,verbatim}

\usepackage[left=1.1in, right=1.1in, top=0.8in, bottom=0.99in]{geometry}

\def\reals{{\mathbb R}}

\def\proof{\noindent{\em Proof: }}
\def\qed{$\spadesuit$}
\def\A{{\cal A}}
\def\B{{\mathcal B}}
\def\F{{\mathcal F}}
\def\G{{\mathcal G}}

\begin{document}
\begin{titlepage}
\title{A Crossing Lemma for Jordan Curves\footnote{Results of this paper have been partly reported in the Proceedings of the 27th Annual ACM-SIAM Symposium on Discrete Algorithms, \cite{PRT16}.}}

\author{J\'anos Pach\thanks{EPFL, Lausanne and R\'enyi Institute,
    Budapest. Supported by Swiss National
    Science Foundation Grants 200020-162884 and 200021-175977. Email: {\tt pach@cims.nyu.edu}}
  \and Natan Rubin\thanks{Ben Gurion University of the Negev,  Beer-Sheba, Israel.
Email: {\tt rubinnat.ac@gmail.com}. Ben Gurion University of the Negev, Beer-Sheba, Israel. Email: rubinnat.ac@gmail.com. Ralph Selig Career Development Chair in Information Theory. Supported in part
by grant 1452/15 from Israel Science Foundation and by grant 2014384 from the U.S.-Israeli Binational Science Foundation. The project leading to this application has received funding from European Research Council (ERC)
under the European Unions Horizon 2020 research and innovation programme under grant agreement No. 678765.}
  \and G\'abor Tardos\thanks{R\'enyi Institute, Budapest. Supported by the ``Lend\"ulet'' Project in Cryptography of the Hungarian Academy of Sciences and the National Research, Development and Innovation Office --- NKFIH projects K-116769 and SNN-117879. Email: {\tt tardos@renyi.hu}}}

\maketitle
\begin{abstract}
If two Jordan curves in the plane have precisely one point in common, and there they do not properly cross, then the common point is called a {\em touching point}.  The main result of this paper is a Crossing Lemma for simple curves: Let $X$ and $T$ stand for the sets of intersection points and touching points, respectively, in a family of $n$ simple curves in the plane, no three of which pass through the same point. If $|T|>cn$, for some fixed constant $c>0$, then we prove that $|X|=\Omega(|T|(\log\log(|T|/n))^{1/504})$. In particular, if $|T|/n\rightarrow\infty$, then the number of intersection points is much larger than the number of touching points.
\smallskip

As a corollary, we confirm the following long-standing conjecture of Richter and Thomassen: The total number of intersection points between $n$ pairwise intersecting simple closed (i.e., Jordan) curves in the plane, no three of which pass through the same point, is at least $(1-o(1))n^2$.
\end{abstract}

\end{titlepage}

\section{Introduction}

\subsection{Preliminaries}
\noindent{\bf Arrangements of curves and surfaces.} It was a fruitful and surprising discovery made in the 1980s that the Piano Mover's Problem and many other algorithmic and optimization questions in motion planning, ray shooting, computer graphics etc., boil down to computing certain elementary substructures (e.g., cells, envelopes, $k$-levels, or zones) in arrangements of curves in the plane and surfaces in higher dimensions~\cite{Ed87,KLPS86,PaS09,ShA95}. Hence, the performance of the most efficient algorithms for the solution of such problems is typically determined by the {\em combinatorial complexity} of a single cell or a collection of several cells in the underlying arrangement, that is, the total number of their faces of all dimensions.

The study of arrangements has brought about a renaissance of Erd\H os-type combinatorial geometry.
For instance, in the plane, Erd\H os's famous question~\cite{Er46} on the maximum number of times the unit distance can occur among $n$ points in the plane can be generalized as follows~\cite{CEGSW90}: {\it What is the maximum total number of sides of $n$ cells in an arrangement of $n$ unit circles in the plane?} In the limiting case, when $k$ circles pass through the same point $p$ (which is, therefore, at unit distance from $k$ circle centers), $p$ can be regarded as a degenerate cell with $k$ sides.

Several beautiful paradigms have emerged as a result of this interplay between combinatorial and computational geometry, from the random sampling argument of Clarkson and Shor~\cite{CS89} through epsilon-nets (Haussler-Welzl~\cite{HW87}) to the discrepancy method (Chazelle~\cite{Cha00}).
It is worth noting that most of these tools are restricted to families of curves and surfaces of {\it bounded description complexity}. This roughly means that a curve in the family can be given by a bounded number of reals (like the coefficients of a bounded degree polynomial). For the exact definition, see \cite{ShA95}.

\smallskip

Another tool that proved to be applicable to Erd\H os's questions on repeated distances is the {\em Crossing Lemma} of Ajtai, Chv\'atal, Newborn, Szemer\'edi and Leighton~\cite{ACNS,L}. It states that no matter how we a draw a sufficiently dense graph $G=(V,E)$ in the plane or on a fixed surface, the number of crossings between its edges is at least
$
\Omega(|E|^3/|V|^2).
$

In particular, this implies that if $G$ has a lot more edges than vertices (that is, $|E|/n\to\infty$), then its number of crossings is much larger than its number of edges. The best known upper bound on the $k$-set problem~\cite{De98}, needed for the analysis of many important geometric algorithms, and the most elegant proofs of the Szemer\'edi-Trotter theorem~\cite{SzT83a}, \cite{SzT83b} on the maximum number of incidences between a set of points and a set of lines (or other, more complicated, curves) were also established using the Crossing Lemma~\cite{PaS98}. These proofs easily generalize from lines to pseudo-segments (i.e., curves with at most one intersection per pair).

\medskip
\noindent{\bf Tangencies and lenses.}
The Circle Packing Theorem of Koebe, Andreev and Thurston \cite{Koe36,An70,Thu97} implies that any $n$-vertex planar graph is isomorphic to a graph whose vertices correspond to $n$ non-overlapping disks in the plane, two vertices being connected by an edge if and only if the boundary circles of the corresponding disks touch each other. Conversely, for any set of closed Jordan curves in general position in the plane, with the property that any two curves are either disjoint or touch at a single point, the corresponding touching graph is easily seen to be planar (see, e.g., \cite{KLPS86}).

The present work furthers the above relation by showing, in analogy to the Crossing Lemma, that the number of proper crossing points among
$n$ Jordan curves in general position grows faster than the number $\tau$ of touching pairs, provided that $\tau/n\to\infty$.

Previously, the study of tangencies in arrangements of curves has been mostly restricted to special families of curves (e.g., boundaries of convex sets or curves of bounded description complexity). Motivated by potential applications to motion planning, Tamaki and Tokuyama~\cite{TT98} extended the $k$-set bounds and incidence bounds from lines to more general curves, by trying to cut the curves into as few pseudo-segments as possible, and then applying the known bounds to them. In this context, the number of {\em tangencies (touchings)} between the original curves plays a special role. By locally perturbing two curves in a small neighborhood of their touching point, one can create two nearby crossings and a small ``lens'' between them. In order to decompose the curves into pseudo-segments, we have to make at least one cut on the boundary of each lens. In many scenarios, the number of cuts needed is roughly proportional to the number of touching points, more precisely, to the maximum number of non-overlapping lenses. This approach was later refined and extended in a series of papers~\cite{ArS02}, \cite{AgS05}, \cite{Chan1}, \cite{Chan2}, \cite{Chan3}, \cite{MaT06} and \cite{PseudoCircles}.

In particular,  Agarwal {\it et al.}\ \cite{PseudoCircles} studied arrangements of pseudo-discs (that is, closed Jordan curves with at most two intersections per pair) and used lenses to establish several fundamental results on geometric incidences and
cell complexity.
Their analysis crucially relied on the following claim: Any family of $n$ pairwise intersecting pseudo-circles admits at most $O(n)$ tangencies. In the special case where the curves are algebraic, any incidence or tangency can be described by a polynomial equation.
Following the pioneering work of Dvir~\cite{Dv10}, Guth and Katz~\cite{GK10}, \cite{GK15}, many of these problems have been revisited from an algebraic perspective.

The structure of tangencies between {\em convex sets} was addressed in \cite{PST12}. It was shown that the number of tangencies between $n$ members of any family of plane convex sets that can be obtained as the union of $k$ packings (systems of disjoints sets) is at most $O(kn)$. The proof of this fact is somewhat delicate, because the boundaries of two convex sets can cross any number of times.

\medskip
\noindent{\bf Richter-Thomassen Conjecture.} Richter and Thomassen conjectured in 1995~\cite{RiT95} that the total number of intersection points between $n$ pairwise intersecting closed Jordan curves in general position in the plane is at least $(1-o(1))n^2$.

Note that if there are no tangencies between the curves, then any two curves intersect at least twice, so that the number of intersection points is at least $2{n\choose 2}=(1-o(1))n^2$. However, if touchings are allowed, the situation is more complicated.

The best known general lower bound is due to Dhruv Mubayi~\cite{Mu02}, who showed that the number of intersection points is at least $(4/5-o(1))n^2$. If any pair of curves have at most a bounded number of points, then the conjecture follows from the K\H{o}v\'{a}ri--S\'{o}s--Tur\'{a}n Theorem \cite{KST54} in extremal graph theory, as proved by Salazar \cite{Sa99}. In an earlier paper~\cite{PRT15}, the authors settled the special case where the curves are {\em convex} or, more generally, if each curve can be cut into a constant number of $x$-monotone arcs. (An arc is called {\em $x$-monotone} if every vertical line intersects it in at most one point.) The problem has remained open for general families of simple closed curves.

\medskip
\noindent{\bf Algebraic techniques.}
As mentioned before, the polynomial technique of Guth and Katz \cite{GK10,G2}, which led to a spectacular breakthrough concerning Erd\H os's problem on distinct distances, has inspired a lot of recent research related to incidences between points, curves, and surfaces \cite{G3,G1,GK15,SSZ15}.
For instance, Ellenberg, Solymosi and Zahl \cite{ESZ16} have shown that any family of low-degree algebraic curves in the plane determines $O(n^{3/2})$ tangencies (where the constant of proportionality can depend on the maximum degree of the curves). Unfortunately, the new techniques only apply in an algebraic framework, where the curves and surfaces in question must be algebraic varieties of bounded degree. Since two algebraic curves of bounded degree that do not share a component have only a bounded number of points in common, restricting the Richter-Thomassen conjecture to such curves, reduces the question to the above mentioned result of Salazar \cite{Sa99}. For many similar problems related to intersection patterns of curves, including the special case of the Erd\H os-Hajnal conjecture\footnote{The conjecture, as applied to string graphs, claims that there exists $\epsilon>0$ such that among $n$ Jordan curves in the plane one always finds at least $n^\epsilon$ pairwise intersecting or at least $n^\epsilon$ pairwise disjoint ones.} \cite{EH89}, our present techniques are not sufficient to handle the case when two curves may intersect an arbitrary number of times \cite{FPT11,FP08,FP10}. There are only very few exceptional examples, when one is able to drop this assumption \cite{Ma14,FP12}. The main result of this paper represents one of the rare exceptions.

\subsection{Our results} The main result of this paper is a Crossing Lemma for a family of Jordan curves. We are going to show, roughly speaking, that the number of {\em proper} (i.e., transversal) crossings between the curves is much larger than the number of touching pairs of curves, provided that the number of touching pairs is super-linear in the number of curves.

To formulate this result more conveniently, we need to agree on the terminology. We say that two (open or closed) curves {\em intersect} if they have at least one point in common. An intersection point $p$ is called a {\em touching point} (in short, a {\em touching}) if $p$ is the {\em only} intersection point of the two curves, and they do not properly cross at $p$.  Note that this definition is somewhat counterintuitive: we do not call a point of tangency between two curves a touching if the curves also intersect at another point. Without this restriction we cannot claim that there are much more crossings than touching points. Indeed, consider $n$ lines in general position in the plane. Notice that one can slightly perturb them to turn each crossing into a proper crossing and a separate point of tangency. In such an arrangement,  half of the $2{n\choose 2}$ intersection points are tangencies. It is assumed throughout that all curves are in {\em general position}, that is, no three of them pass through the same point and no two share infinitely many points.

We state our Crossing Lemma in two forms. First, we formulate it for pairwise intersecting closed curves. In this formulation, we can prove a slightly better asymptotic gap between the number of intersections and the number of touchings:

\begin{theorem}\label{main}
Let $\A$ be a collection of $n$ pairwise intersecting closed Jordan curves in general position in the plane. Let $T$ denote the set of touching points and let $X$ denote the set of intersection points between the elements of $\A$. We have\footnote{The preliminary version \cite{PRT16} states a too optimistic (and, unfortunately, less accurate) estimate of $\Omega\left((\log\log n)^{1/8}\right)$.}
\begin{equation*}
\frac{|X|}{|T|}=\Omega\left((\log\log n)^{1/12}\right).
\end{equation*}
\end{theorem}

\medskip
We will see that Theorem \ref{main} is an easy corollary to its bipartite version:

\begin{theorem}\label{main-bipartite}
Let $\F$ and $\G$ be two disjoint collections of closed Jordan curves in the plane, each consisting of $n$ curves. Suppose that $\F\cup\G$ is in general position and that any two curves from the same collection intersect. Let $T$ denote the set of all touching points between the curves in $\F$ and the curves in $\G$, and let $X$ denote the set of all intersection points between the elements of $\F\cup \G$. We have
\begin{equation*}
\frac{|X|}{|T|}=\Omega\left((\log\log n)^{1/12}\right).
\end{equation*}
\end{theorem}

Most of this paper is devoted to the proof of Theorem~\ref{main-bipartite}. Theorem~\ref{main} can be deduced from Theorem~\ref{main-bipartite}, as follows.
\medskip

\noindent{\em Proof of Theorem \ref{main} (using Theorem~\ref{main-bipartite}):}
Assume without loss of generality that $n$ is even and consider a random partition of the curves of $\A$ into $n/2$-sized families $\F$ and $\G$. Notice that the expected number of the original touchings with the two touching curves ending up in distinct families $\F$ and $\G$ is at least $|T|/2$.
The overall number of intersection points does not change, so applying Theorem~\ref{main-bipartite}) to the families $\F$ and $\G$ proves the statement of Theorem~\ref{main}.
\qed

\medskip
We use Theorem~\ref{main} to settle the Richter-Thomassen conjecture~\cite{RiT95}:

\begin{theorem}\label{richterthomassen}
The total number of intersection points between $n$ pairwise intersecting closed curves in general position in the plane is at least $(1-o(1))n^2$.
\end{theorem}

\medskip

\noindent{\em Proof of Theorem~\ref{richterthomassen} (using Theorem~\ref{main}):}  It is enough to notice that if $|T|=o(n^2)$, then the statement follows from the trivial bound $|X|\ge2{n\choose2}-|T|$. Otherwise, if $|T|\ge \varepsilon n^2$ for some $\varepsilon>0$, Theorem~\ref{main} immediately implies that $|X|=\Omega(n^2(\log\log n)^{1/12})$, which is much better than required. \qed

\medskip
The assumption that the curves are closed is crucial for Theorem \ref{richterthomassen}, as a family of pairwise-intersecting segments in general position in the plane determines only ${n \choose 2}<n^2/2$ intersections.
However, both Theorem \ref{main} and Theorem \ref{main-bipartite} readily extend to families of Jordan arcs. Indeed, slightly inflating the Jordan arcs to closed Jordan curves, while preserving all of the touching pairs, we increase the number of intersection points by a factor of at most $4$.

\medskip
Next, we formulate a version of Theorem~\ref{main} without the assumption that the curves are pairwise intersecting. Note, however, that one may draw $n$ circles with up to $3n-6$ touchings and no proper crossings~\cite{PseudoCircles,PaS09}.
One might believe that some linear lower bound on the number of touching, like $|T|\ge10 n$ should be enough for us to prove a separation $|T|=o(|X|)$, but this is false as shown by the following example. Fix a large constant $k$ and consider $n-k$ pairwise disjoint unit circles in the plane.
It is easy to select $k$ other closed curves in general position such that each of them touches every circle and any pair of them intersect at most $n-k+1$ times. In this arrangement, $|T|=k(n-k)$ and $|X|\le|T|+{k\choose 2}(n-k+1)$, so that we have $\frac{|X|}{|T|}\le k$, a constant. This motivates to formulate our lower bound on $\frac{|X|}{|T|}$ not as a function of $n$ assuming a lower bound on $|T|$, but rather as a function of $|T|/n$. We conjecture the following strong bound. We formulate it for Jordan arcs and not for closed Jordan curves.
\medskip

\begin{conjecture}\label{conjecture} Let $\A$ be a collection of $n$ Jordan arcs in general position in the plane. Let $T$ denote the set of touching points and $X$ the set of intersection points between the elements of $\A$. We have
$$\frac{|X|}{|T|}=\Omega\left(\log\frac{|T|}n\right).$$
\end{conjecture}

The conjectured logarithmic separation between $X$ and $T$, if true, cannot be improved. Indeed, Fox {\em et al.}\ \cite{FFPP10} constructed two $n$-sized families, $\F$ and $\G$, of pairwise intersecting $x$-monotone curves in the plane such that every curve in $\F$ touches every curve in $\G$, and the total number of intersections between the members of $\F\cup \G$ is $O(n^2\log n)$. They also showed that for this setting one always has $\Omega(n^2\log n)$ intersections.

Though we have been unable to verify Conjecture~\ref{conjecture}, we can deduce the following weaker bound from Theorem~\ref{main-bipartite}.

\begin{theorem}\label{general}Let $\A$ be a collection of $n$ Jordan arcs in general position in the plane. Let $T$ denote the set of touching points and $X$ the set of intersection points between the elements of $\A$. We have
$$\frac{|X|}{|T|}=\Omega\left(\left(\log\log\frac{|T|}n\right)^{1/504}\right).$$
\end{theorem}

\subsection{Organization and overview}
The paper is organized as follows.

In Section \ref{Sec:MainProof}, we prove our main technical tool, Theorem~\ref{main-bipartite}. Unlike many previous bounds on the crossing numbers of geometric structures, which relied on Euler's formula \cite{ACNS,L} or parity arguments from topology \cite{Tut70},  our analysis is based on a more local machinery of {\it charging schemes} -- a powerful yet simple method developed in Computational Geometry to estimate the number of special features of bounded description complexity in arrangements of algebraic curves in $\reals^2$ and surfaces in $\reals^d$; see \cite[Section 7]{ShA95} for a comprehensive demonstration of this technique.

In the most typical planar scenario, we are given an arrangement of $n$ algebraic curves and seek a non-trivial upper bound on the number of ``special" vertices which satisfy a certain topological condition (e.g., vertices that lie on the boundary of a given face).
That is, we are to show that the concerned vertices are relatively scarce, and the vast majority of the intersection points do not possess the desired property.
To this end, we assign each special vertex $v$ to several other vertices $v'$ in the arrangement. The assignment is fractional and specified by a rule in which $v$ ``receives" at least $c_{\mbox{\scriptsize in}}$ units of charge from the other vertices $v'$.
In most instances, the charging rule is of a fairly local nature and respects some natural criterion of proximity between $v$ and $v'$ within the arrangement. A successful charging scheme must guarantee that the total charge ``sent" by any vertex $v'$ is much smaller than $c_{\mbox{\scriptsize in}}$.

We adapt the above charging paradigm to show that  only few of the intersection points can be touchings. Since the original machinery applies (with very few exceptions) only to objects of bounded description complexity, the adaptation requires extra care
to avoid ``overcharging" of intersection points.


In Section \ref{Sec:Dense}, we establish Theorem~\ref{general}.
The proof proceeds as follows.
First, we use a simple sampling argument to
replace the hypothesis of Theorem \ref{main} with a somewhat weaker one -- the intersection graph must be rather dense.
To deduce Theorem \ref{general} from the ``dense" Crossing Lemma, we partition the arrangement into sufficiently dense pieces by repeatedly applying the separator result of Fox and Pach \cite{FP08}.

\section{Proof of Theorem \ref{main-bipartite}}\label{Sec:MainProof}

This is the most complex part of the paper. We start with a brief and informal outline of the proof. First we bound the number of touching points $t\in T$ that are contained in an arc of arbitrary size whose ``crossing to touching ratio'' is high. We will call these ``happy'' touching points and bounding their number relative to the number $|X|$ of intersection points is simple. The rest of the touching points we call ``sad'' and denote their set by $T'$. We use the so-called charging method to bound $|T'|$: we send certain amounts of ``charge'' from points in $X$ to points in $|T'|$. If we
manage to make sure that the total charge sent by any point in $X$ is at most
$c_{\mbox{\scriptsize out}}$ and the total charge received by any point $t\in T'$ is
at least $c_{\mbox{\scriptsize in}}$, then we have established that $|X|/|T'|\ge
c_{\mbox{\scriptsize in}}/c_{\mbox{\scriptsize out}}$. Note that we used the same method in our paper \cite{PRT15}
to prove certain special cases of the Richter--Thomassen conjecture. For most of the charging argument, there is no need to restrict our attention to sad touching points, but at one crucial point, namely in the proof of Claim~\ref{close}, it helps us that we have already taken care of all happy ones.
\medskip

The proof of Theorem~\ref{main-bipartite} is organized as follows.
In Section~\ref{subsec:happy}, we fix the arrangement of curves, define happy and sad touching points and bound the number of happy touching points relative to the number of intersection points.

In Section~\ref{subsec:rules}, we define our three charging rules by which intersection points send specified amount of ``charge" to sad touching points. These rules involve a parameter $k$ that we call the ``scale''. The charging rules should be applied in many phases, each phase with a different scale.

In Section~\ref{Subsec:UpperBound}, we specify the number $M$ of phases and the corresponding values of the scale parameter $k$. We establish a constant upper bound on the amount of charge sent by any point of $X$, averaged over the $M$ phases. For the first and third charging rules, we have a constant upper bound in each individual phase, but for the second rule this is not the case, and here the appropriate choice of scales is important. In the rest of the proof, it is almost irrelevant how we set the value of the scale parameter $k$ (within reasonable limits).

Section~\ref{Subsec:LowerBound} is devoted to proving a lower bound of $\Omega((\log\log n)^{1/12})$
on the amount of charge received by a sad touching point in any given phase.
We will use Claim~\ref{close} in this proof, whose proof is technically
involved and is postponed to Section~\ref{Subsec:ProofClose}.

In Section~\ref{wrap}, we finish the proof of Theorem~\ref{main-bipartite} by bounding the number of sad touching points.

\subsection{Happy and sad touchings}\label{subsec:happy}

Let us fix the families $\F$ and $\G$ of $n$ closed Jordan curves, as in the statement of Theorem~\ref{main-bipartite}. For the sake of brevity we write $\A=\F\cup \G$.
Let $\alpha_1=(\log\log n)^{1/12}/10$, where $\log$ denotes the binary
logarithm.
We need to prove $|X|/|T|=\Omega(\alpha_1)$. For this proof we
assume that $|X|\le \alpha_1n^2$ as otherwise the statement follows from the trivial bound $|T|\le n^2$. Note that as the statement we want to prove is asymptotic we
can simplify our calculations by always assuming that $n$ is large enough.

We call an arc $a^*$ contained in one of the curves $a\in\A$ {\em happy} if $|X\cap a^*|\geq
\alpha_1|T\cap a^*|$. We say that a touching point $t\in T$  is {\em happy} if  it is contained in a happy arc, otherwise we call $t$ {\em sad} and denote the set of sad touching points by $T'$.

\begin{figure}[htbp]
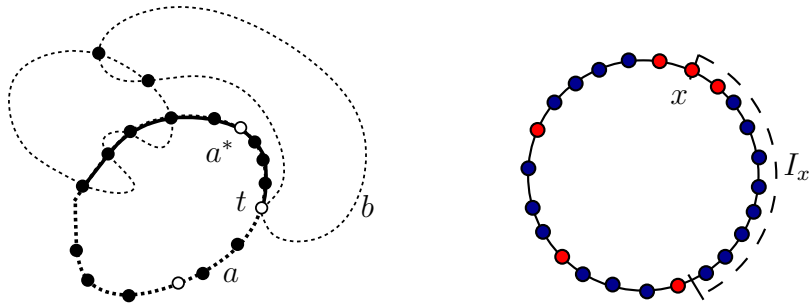

\begin{center}
\input{HappyTouching.pstex_t}\hspace{2cm}\input{Covering.pstex_t}
\caption{\small Left -- The touching $t$ is happy. The hollow points belong to $T\cap a$. The arc $a^*\subset a$
  satisfies $|X\cap a^*|\geq \alpha_1 |T\cap a^*|$ with $\alpha_1=5$. Right -- Lemma
  \ref{Lemma:RedBlue}. Each red point $x\in R$ is contained in an arc $I_x$
  that satisfies $w(B\cap I_x)\geq \lambda |R\cap I_x|$. (In the depicted
  scenario $w(x)=1$ for all $x$ and we have $\lambda=4$.)}
\label{Fig:Covering}
\end{center}
\end{figure}

To bound the number of happy touching points, we use the following simple lemma. See Figure \ref{Fig:Covering} (right). Notice that the sets $R$ and $B$ may overlap. In the use of this lemma below the weight function $w$ is constant. We still formulate the lemma with an arbitrary positive weight function $w$ because we will use the same lemma a second time later in this section and there we will use a non-constant weight function.

\begin{lemma}\label{Lemma:RedBlue}
Let $a$ be a simple (open or closed) Jordan curve and let $R$ and $B$ be two finite subsets of points on $a$. Let $\lambda$ be a positive
constant and let $w: B\rightarrow \reals$ be a
positive weight function. For $S\subseteq B$ the weight of $S$ is $w(S)=\sum_{x\in S}w(x)$. If every point $x\in R$ is contained in an arc $I_x\subseteq a$ that satisfies $w(B\cap I_x)\geq \lambda |R\cap I_x|$,
then we have $w(B)\geq \lambda|R|/3$.
\end{lemma}

\medskip
\proof We prove by induction on $|R|$. The claim trivially holds if $R$
is empty, so we assume $R$ is not empty and the statement of the lemma holds
for $R'$ and $B'$ as long as $|R'|<|R|$.

Let us choose $x\in R$ to maximize $|R\cap I_x|$ breaking ties arbitrarily. Let
$B'=B\setminus I_x$ and $R'=\{y\in R\mid I_y\cap I_x=\emptyset\}$. For $y\in
R'$ we have $B'\cap I_y=B\cap I_y$ and $R'\cap I_y\subseteq R\cap I_y$, so the
assumption of the lemma is satisfied for $R'$ and $B'$. As $x\notin R'$ we have
$|R'|<|R|$ and thus, by the inductive hypothesis, we have
$w(B')\ge\lambda|R'|/3$. By the choice of $x$ every  $y\in R\setminus R'$
must either be in $I_x$ or it is one of the $|I_x\cap R|$ next points in $R$
in either side of the arc $I_x$. So we have $|R|-|R'|\le3|I_x\cap R|$. We
further have $w(B)-w(B')=w(I_x\cap B)\ge\lambda|I_x\cap
R|\ge\lambda(|R|-|R'|)/3$. Adding this inequality to the one obtained from the
inductive hypothesis finishes the proof. \qed
\medskip

\begin{lemma}\label{happy}
Let $T$ and $X$ be the respective sets of touching points and intersection points as defined in Theorem \ref{main-bipartite}, and let $T'\subset T$ be the set of sad touching points. Then we have $|T|-|T'|\le 6|X|/\alpha_1$.
\end{lemma}
\medskip

\proof We apply Lemma \ref{Lemma:RedBlue} for each curve $a\in \A$ with
$\lambda=\alpha_1$, $B=B_a=a\cap X$ and $R=R_a$ being the set of touching points contained in a happy arc $a^*\subset a$. We use the uniform weight function $w(x)=1$ for each $x\in
B_a$. We obtain $|B_a|\ge\alpha_1|R_a|/3$. Summing this for all
$a\in \A$ we get $2|X|$ on the left hand side and at least $\alpha_1(|T|-|T'|)/3$
on the right hand side. \qed

\subsection{The charging rules}\label{subsec:rules}

As mentioned in the outline above we bound $|T'|$ using a charging scheme.

Each charging scheme describes a fractional assignment of the elements of a set $A$ to the elements of another set $B$, and can be described as a weight assignment to the edges of the complete directed bipartite graph $B\times A$.
In the language of charging schemes, $a\in A$ {\it receives} $w(b,a)$ units of charge from $b$, whilst $b$ {\it sends} $w(a,b)$ units to $a$.

The eventual upper bound on $|A|$ in terms of $|B|$ depends on the minimal weighted indegree $c_{\mbox{\scriptsize in}}=\min_{a\in A}\sum_{b\in B} w(b,a)$ and the maximum weighted outdegree $c_{\mbox{\scriptsize out}}=\max_{b\in B}\sum_{a\in A}w(b,a)$. With these parameters, a standard double counting argument shows that $\frac{|A|}{|B|}\le\frac{c_{\mbox{\scriptsize out}}}{c_{\mbox{\scriptsize in}}}$.

Our charging is done in phases: in each phase we fix the value of the parameter $k$ (``the
scale'') and perform certain chargings with that scale. Our goal is to make sure that each point
in $X$ sends out a constant charge in each phase, while each touching in $T'$
receives a charge of $\Omega(\alpha_1)$ in each phase. If we could do this, then a single phase would be enough
to prove Theorem~\ref{main-bipartite}. But we will not quite achieve this goal.
Some points in $X$ will be overcharged in certain
phases: they send out more than a constant amount of charge. This problem is solved by considering
several phases at once.
The exact values of the scale parameter $k$ will be set in Section~\ref{Subsec:UpperBound} to ensure that {\em on average} no intersection is overcharged.

\paragraph{Arcs and lenses.} Before specifying the exact charging rules we introduce some notation. We orient each curve $a$ in $\F$ so that all other curves from $\G$ touching $a$ touches it on its right side. This
is possible as if $a\in \F$ has a touching curve on either side, then these curves
are not intersecting counter to our assumption that each pair of curves in
$\G$ intersect. We similarly orient the curves of $\G$. We use the word {\em arc} for closed segments of the curves in $\A$. We will use lowercase
letters with an asterisk to denote arcs. The arcs inherit
their orientation from the curve of $\A$ containing them and this orientation distinguishes the
{\em starting point} and the {\em end point} of an arc. For distinct points $p$ and $q$ in a curve
$a\in \A$ we write {\em the arc of $a$ from $p$ to $q$} to refer to the single arc on $a$ with $p$
as its starting and $q$ as its end point. We can simply refer to an arc as ``the arc from $p$ to $q$''
unless $p,q\in X$ represent two intersections of the same two curves from $\A$.
The orientation also makes references like ``the next $k$
points of $T$ along $a$ after $p$'', or ``the last $k$ points of $T$ along $a$ before $p$''
unambiguous. By the {\em length} of an arc
we mean the number of sad touching points it contains. Let $x\in X$ be a non-touching intersection point of the curves $a,b\in\G$, and let $y$ be another intersection point of the same two curves, the next such point along
$a$.

We call the arc of $a$ from $x$ to $y$ a {\em lens}. (Note that most texts include
both arcs from $x$ to $y$ in their definition of a lens, but for us it is
simpler to focus on a single arc. We will only use the term lens for lenses determined by curves in $\G$.)

\medskip
We set the following parameters: $\alpha=\alpha_1+2$, $v=21000\alpha^{12}$
depending only on $n$ and the parameter $w=w(k)=k^3/(4000\alpha^5n^2)$ that
also depends on the scale.

\begin{figure}[htbp]
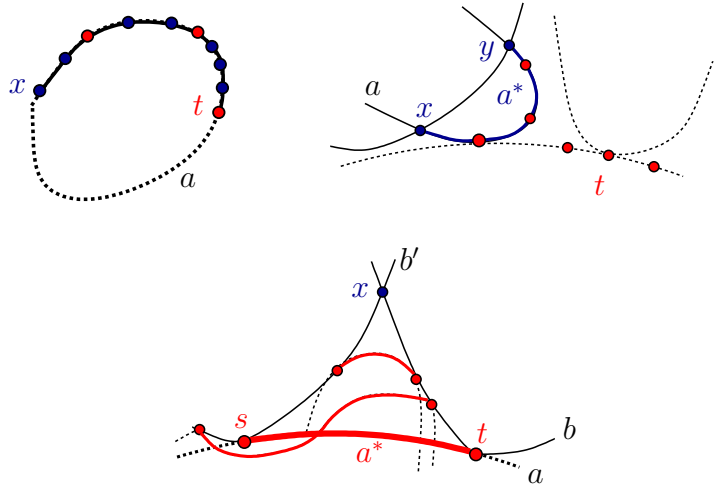

\begin{center}
\input{FirstRule.pstex_t}\hspace{1.3cm}\input{ChargingRuleLens.pstex_t}\\

\vspace{0.5cm}\input{LastRuleApex.pstex_t}
\caption{\small Left: First charging rule. A point $x\in X$ sends $1/k$ units
  to a sad touching $t\in T'$ if the interval from $t$ to $x$ (or vice versa)
  has length at most $k$. Right: Second charging rule. The sad touching $t$
  receives $v/(k(l+w))$ units from the lens $a^*$ with endpoints $x$ and
  $y$. Bottom: Third charging rule. The starting point $t$ of the poor arc $a^*$
  receives $2a/k$ units from the apex $x\in X$ because there exist at most $k/a$ other poor arcs with apex $x$ that start between $x$ and $t$.}
\label{Fig:Rules}
\end{center}

\end{figure}

Let us consider the phase with scale $k$. We start with describing our three
charging rules sending charges from intersection points in $X$ to sad
touchings in $T'$. See Figure \ref{Fig:Rules}.

\medskip
\noindent {\bf First charging rule.} A point $x\in  X$ sends a charge of $1/k$ to a
point $t\in T'$ if the two points are on a common curve of $\A$ and
either the arc from $x$ to $t$ or the arc from $t$ to $x$ has length at most $k$.

\medskip
\noindent{\bf Second charging rule.} If the length $l$ of a lens $a^*$ satisfies
$l\le3\alpha^3k$, then $a^*$ sends a charge of $v/(k(l+w))$ to all
points
$t\in T'$ that have an arc of length at most $k+1$ along a curve in $\F$ from $t$ to a point in
$a^*\cap T'$.

\smallskip
For accounting purposes, we consider a charge sent by a lens $a^*$ to be sent
by the starting point of $a^*$. Note that exactly two lenses starts at every
non-touching intersection point between two curves from $\G$.

We call a point of $T'$ {\em poor} in this phase if it receives less than a
total charge of $\alpha$ from the first two charging rules. We call an arc
{\em poor} if it starts at a poor touching point,
ends at a sad touching point and has length at most $k+1$.

Let $a^*$ be an arc of a curve $a\in
\F$ starting at $t\in T'$ and ending at $s\in T'$. Let $b$ and $b'$ be the curves in $\G$ touching
$a$ in the points $t$ and $s$, respectively. We define the {\em apex} of the arc
$a^*$ as the first point on $b'$ after $s$ that also belongs to $b$. This is a
well defined point in $X$ as $b$ and $b'$ (as any pair of curves in $\G$) must
intersect.

\medskip
\noindent{\bf Third charging rule.} Let $a^*$ be a poor arc starting at $t\in T'$ and
having $x\in X$ as its apex. The intersection point $x$ sends a charge of $2\alpha/k$ to $t$ in this phase unless there are
more than $k/\alpha$ poor arcs, each starting at a point in the arc from $x$
to $t$ and having $x$ as its apex.

\subsection{Total charge sent}\label{Subsec:UpperBound}

\begin{lemma}\label{out13}
The total charge sent from a intersection point $x\in X$ in a phase according
to the first and third rules is at most $8$.
\end{lemma}

\proof The first rule sends a charge of $1/k$ from $x$ to the first $k$ sad
touching points in each of four ``directions'' (in both directions of both
curves containing $x$). That is at most $4k$ sad touching points
for a total charge of at most $4$.

The third rule sends a charge of $2\alpha/k$ to the first $\lfloor k/\alpha\rfloor$ touching points from $x$
along either curves containing $x$ satisfying a certain condition (namely being the starting point of
a poor arc having $x$ as its apex) for a total charge of at most $4$.
\qed

\smallskip
Note that both the first and the third rule charges an intersection point of two curves in $\G$ irrespective of the scale. In contrast, each lens has an intrinsic length value $l$, which roughly describes the scale $k$ of the phase where this vertex can be charged via the second rule.
A statement similar to Lemma~\ref{out13} is false for the second charging rule because it severely overcharges
the lenses whose length is approximately $k$. Observe, however, that the rule does not
charge a lens longer than $3\alpha^3k$ and charges it very lightly if the lens is much
shorter than $w=w(k)$. This is enough for us to set up the scales of the different phases in such a way
that no intersection point is overcharged on average.

For technical reasons, in phase $k$ we overcharge lenses of length $l$ within a fairly large interval, namely for $k^3/(n^2 {\it poly}(\alpha))<l<k \,{\it poly}(\alpha)$. To avoid overcharging the same lens in many phases, we can only have $O(\log\log n)$ phases. This is the reason that our lower bound on $|X|/|T|$ in Theorem~\ref{main-bipartite} (and as a consequence also in Theorems~\ref{main} and \ref{general}) is substantially weaker than the similar bound in \cite{PRT15}.

We use the following scales for the different phases of our charging: $k=80\alpha^4 n/2^{3^u}$, where $u$ is an
integer satisfying $(\log\log n)/5<u\le (\log\log n)/2$. We have $M=\lfloor(\log\log n)/2\rfloor-\lfloor (\log\log
n)/5\rfloor$ phases.

\begin{lemma}\label{out2} For any intersection point $x\in X$, the charge leaving $x$ by the second rule averaged over the $M$ phases is at most $2$.
\end{lemma}

\proof Each non-touching intersection point of two curves of $\G$ is the starting point of at most two
lenses (no lens starts at a touching point). We bound the average charge sent by a fixed lens $a^*$ by $1$. Let
$l$ be the length of $a^*$ and let $k_0$ be the smallest scale of a phase
where the lens $a^*$ is charged. We have
$l\le3\alpha^3k_0$ and the total charge $a^*$ sends in this phase is less than $v$. For phases with scale
$k>k_0$, we bound the charge $a^*$ sends by
$vl/w(k)\le 3v\alpha^3k_0/w(k)$.
With the prior choice of scales $k$ and the parameter $w(k)$, we have
$w(k)\geq 3\alpha^3k_0$ and the value of $w(k)$ grows by a factor greater than
$2$ every time we go from a scale to a larger scale. Thus, the total
charge $a^*$ sends in all the phases is at most $3v$. With our choice of the
parameters, we have $M\ge3v$ and this proves the estimate claimed. \qed

\subsection{Total charge received.}\label{Subsec:LowerBound}

Our goal in this section is to prove the following lemma.

\begin{lemma}\label{in}
Every sad touching point receives a total charge of at least $\alpha$ in
every phase.
\end{lemma}

We start with an
informal summary of the argument.

Let $t\in T'$ be a sad touching point between a pair of curves $a\in \F$ and $b\in \G$. We consider the sequence $t_1,\ldots, t_k$ of the first $k$ sad touching points that follow $t$ along $a$.

We can assume that $t$ is poor, as otherwise it receives enough charge by the first
two rules. This implies, in particular, that the arc $a_k^*$ of $a$ from $t$
to $t_k$ (of length exactly $k+1$) contains fewer than $\alpha k$ intersection
points of $X$, as all these points send a charge of $1/k$ to $t$ by the first
rule.

\begin{figure}[htbp]
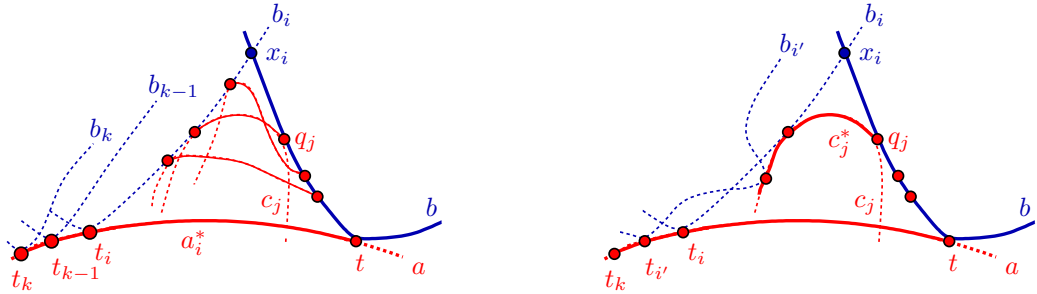

\begin{center}
\input{Poor.pstex_t}\hspace{2cm}\input{Poor1.pstex_t}
\caption{\small Left -- The poor point $t$ receives less than $\alpha$ units of charge through the third charging rule. For most $1\leq i\leq k$ the arc $a_i^*$ shares its apex $x_i$ with $k/\alpha$ poor arcs whose underlying curves $c_j$ meet $a_i^*$. Right -- The curve $c_j\in Q_i\cap Q_{i'}$ is simultaneously tangent to both $b_i$ and $b_{i'}$.}
\label{Fig:Poor}
\end{center}
\end{figure}

For each $1\leq i\leq k$, we consider the apex $x_i$ of the arc $a_i^*$ of length $i+1$ from $t$ to $t_i$; see Figure \ref{Fig:Poor} (left). Notice that $x_i$ is an intersection point of $b$ and another curve $b_i\in \G$ touching $a$ at $t_i$.

We further assume, for contradiction, that $t$ receives less than $\alpha$ units of charge by the third rule from the apex points $x_i$. This implies that at least half of these points $x_i$ do not send charge to $t$ by the third rule. For each of these points $x_i$, not sending charge to $t$,
there exist more than $k/\alpha$ poor arcs with apex $x_i$, each starting at a poor point $q_j$ on the portion of $b$ from $x_i$ to $t$.
We use $Q_i$ to denote the set of such poor points $q_j$ that are associated with $x_i$.
For each of these poor arcs, its underlying curve $c_j$ must intersect $a$ within $a_i^*$ or, else, it would be trapped in the region of $\reals^2\setminus (a\cup b\cup b_i)$ to the right of the arc $a_i^*$ (and, thereby, remain disjoint from $a$).
This implies that the total number of such curves $c_j$ associated with at least one of the apexes $x_i$ cannot exceed $\alpha k$. Since the overall number of touchings between the curves $b_i$ and $c_j$ is at least $(k/2)\cdot(k/\alpha)$, the resulting bipartite graph of tangencies has density at least $1/(2\alpha^2)$.

All of the above touching points between $c_j\in \F$ and the curves $b_i\in
\G$ must lie within an arc $c_j^*$ of $c_j$ of length $k+1$ which starts at
$q_j$. Using that $q_j$ is poor, the arc $c_j^*$ contains at most $\alpha k$
points of $X$.

Since the graph of touchings between the curves $b_i\in \G$ and $c_j\in \F$ is dense, for an average pair $1\le i,i'\le k$
there exist $|Q_i\cap Q_{i'}|=\Omega^*(k)$ curves $c_j$ that are simultaneously tangent to both $b_i$ and $b_{i'}$.\footnote{The $O^*()$ and $\Omega^*()$ notation hides multiplicative factors of $\alpha$. This notation is only used in this informal proof sketch.} See Figure \ref{Fig:Poor} (right).

\begin{figure}[htbp]
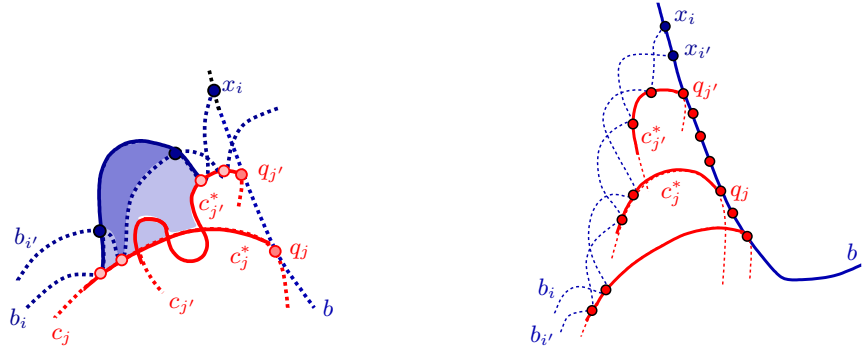

\begin{center}
\input{Lens.pstex_t}\hspace{3cm}\input{Comb.pstex_t}
\caption{\small Left -- In the first charging scenario, the poor touchings $q_j,q_{j'}\in Q_i\cap Q_{i'}$ are close along $b^*$ (so they can meet one another). As $b_{i'}$ enters the pocket formed by the touchings between $b_i$ and each of the curves $c_{j},c_{j'}$, the curves $b_i$ and $b_{i'}$ define a lens of length $l=O^*(k)$ which can send $\Omega^*(1/k^2)$ units to $q_j$. Right -- In the second charging scenario,  the pairwise disjoint arcs $c_j\in Q_i$ form the teeth of a comb-like arrangement.}
\label{Fig:ChargeLens}
\end{center}
\end{figure}

Our parameters are fine-tuned so as to interpolate between the following
extreme scenarios:

\medskip
(i) Any two points $q_j,q_{j'}\in Q_i\cap Q_{i'}$ are close along $b$
in the following sense: the arc of $b$ from $q_j$ to $q_{j'}$, or the complementary arc from $q_{j'}$ to $q_j$, contains at most $2\alpha k$ points of $T$.
Then the two curves $b_i,b_{i'}\in \G$ form a lens of length $l=O^*(k)$ as $b_{i'}$ enters the pocket formed by the touchings between $b_{i}$ and each of the curves $c_j$ and $c_{j'}$, or vice versa. See Figure \ref{Fig:ChargeLens} (left).
The resulting lens of $b_i$ and $b_{i'}$ can send, by our second rule, $\Omega^*(1/k^2)$ units of charge to one of the touchings $q_j,q_{j'}$.
Repeating this argument for $\Omega(k^2)$ pairs
$1\le i,i'\le k$ would eventually contradict the choice of $q_j$ as poor touching points.

(ii) No two points $q_j,q_{j'}\in Q_i\cap Q_{i'}$ are close along $b$.
We argue that, for any $q_j,q_{j'}\in Q_i\cap Q_{i'}$, the respective short arcs $c^*_j\subset c_j$ and $c^*_{j'}\subset c_{j'}$ are disjoint and, therefore, they constitute the ``teeth'' of the comb-like arrangement $\Gamma$ of these arcs together with $b$; see Figure \ref{Fig:ChargeLens} (right).
It then follows that $b_i$ and $b_{i'}$ experience at least $\Omega^*(k)$ intersections, as they touch the neighboring pairs of the teeth of $\Gamma$.
Repeating this for $\Omega(k^2)$ pairs $1\le i,i'\le k$, would contradict the initial assumption that the total number of intersection points satisfies $|X|\leq \alpha_1n^2$.
\medskip

We make the above argument formal and prove Lemma~\ref{in} through a series of small claims. We start with a simple observation that will allow us to speak about ``the
next $k$ sad points'' after a poor point on a curve:

\begin{claim}\label{gek}
If a curve $a\in \A$ contains at most $k$ sad points, then none of them is poor.
\end{claim}

\proof Clearly, if $|a\cap T'|\le k$, then every intersection point in $X\cap
a$ sends a charge of $1/k$ to every point in $T'\cap a$ according to the first
rule. The claim follows as there are at least $n-1$ intersection points on $a$. \qed

\smallskip
For the proof of Lemma~\ref{in} we fix the phase with scale $k$ and we also
fix a single sad  touching point $t\in T'$. We assume for contradiction that
$t$ receives a total charge of less than $\alpha$.
Note first that our assumption implies that $t$ is poor.

Let the curves touching at $t$ be $a\in \F$ and $b\in \G$. Let $t_1,t_2,\dots,t_k$ be the first
$k$ sad touching points after $t$ along $a$. By Claim~\ref{gek}, these exist. For $1\le i\le k$, let $a_i^*$ be the arc of $a$ from $t$ to $t_i$, let
$x_i$ be the apex of $a_i^*$ and let $b_i$ be the curve in $\A$ that touches $a$ at $t_i$. 

\begin{figure}[htbp]
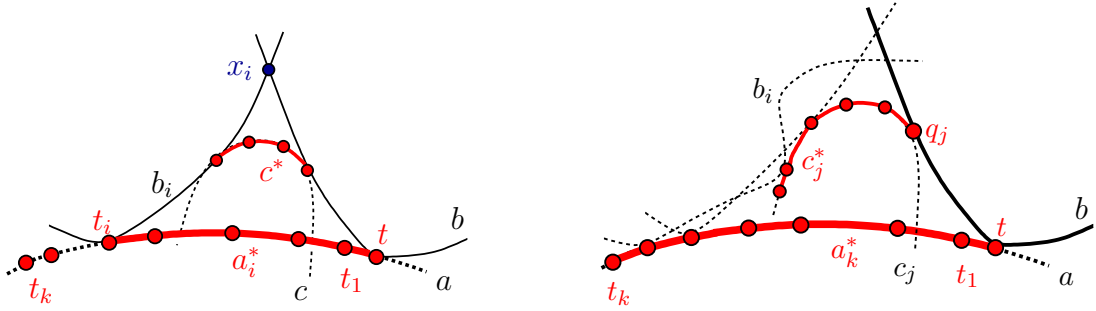

\begin{center}
\input{FastArc.pstex_t}\hspace{1.5cm}
\input{GoodPoint.pstex_t}
\caption{\small Left -- A poor arc $c^*$ is $i$-fast if it has apex $x_i$ and starts on the arc of $b$ from $x_i$ to $t$. Right -- A good point $q_j$ together with the adjacent arc $c_j^*\subset c_j$. Notice that $c_j$ must meet $a$ within $a^*_k$, and $c_j^*$ contains all the fast arcs that start at $q_j$.}
\label{Fig:GoodPoint}
\end{center}
\end{figure}

\medskip
\noindent{\it Defintion.} We call a poor arc
{\em $i$-fast} if it starts at a point in the arc from $x_i$ to $t$ and has $x_i$ as its apex, see Figure \ref{Fig:GoodPoint} (left). We call an arc
{\em fast} if it is $i$-fast for some $1\le i\le k$.

By the third charging rule, each apex point $x_i$ either
transfers $2\alpha/k$ units of charge to $t$ or gives rise to $k/\alpha$ $i$-fast arcs.

\begin{claim}\label{fast} There are more than $k^2/(2\alpha)$ fast arcs.
\end{claim}

\proof Note that $a_i^*$ itself is $i$-fast and, therefore, $t$ receives a charge of $2\alpha/k$ from
$x_i$ according to
the third charging rule, unless there are more than $k/\alpha$ $i$-fast arcs. As $t$ receives a total charge
of less than $\alpha$, we must have more than $k/\alpha$ $i$-fast arcs for each of more than $k/2$
different values of $i$. This proves the claim. \qed
\medskip

Note that all fast arcs start at a sad touching point on $b$, and $k$ of them start at $t$. We call a point
{\em good} if at least $k/(4\alpha^2)$ fast arcs start there. Let us name the good points $q_1,\dots,q_L$
in the order they appear on $b$ starting at $q_1=t$ and going along $b$ in the reverse direction. For
$1\le j\le L$, let $c_j\in \A$ be the curve that touches $b$ at $q_j$ and let $c_j^*$ be the unique arc
from $q_j$ to a point in $c_j\cap T'$ of length exactly $k+1$. The existence
follows from Claim~\ref{gek}. In particular, we have $c_1^*=a_k^*$. See Figure \ref{Fig:GoodPoint} (right).

\begin{claim}\label{goodpoor}
With the previous notation, the following is true.
\begin{enumerate}
\item All good points are poor.

\item We have $|c_j^*\cap X|<\alpha k$ for all $1\le j\le L$.

\item The number of good points is $L\le\alpha k$.

\item At least $k^2/(4\alpha)$ fast arcs start at a good point.

\item Any $i$-fast arc that starts at one of the good points $q_j$ ends at the point where $b_i$ touches
$c_j$, and it is contained in $c_j^*$.
\end{enumerate}
\end{claim}

\proof The first statement holds, because any fast arc starts at a poor point, by definition.

The second statement follows, as each point in $c_j^*\cap X$ sends a charge of
$1/k$ to the poor point $q_j$, by the first rule.

We prove the third statement in a stronger form: the same bound holds for the number $L'$ of all
starting points of fast arcs. Let $c^*$ be an $i$-fast arc starting at $q\ne t$. The curve $c\in \A$ containing
$c^*$ must intersect $a$ and, thus, it must escape the triangle like region
bounded by the arc of $b$
from $x_i$ to $t$, $a_i^*$ and the arc of $b_i$ from $t_i$ to $x_i$. It cannot
cross $b$ or $b_i$, so it must leave through (or touch) $a_i^*$, ``using up''
at least one of the at most $\alpha k$ intersection points on
$a_k^*$, as $a_i^*$ is contained in $a_k^*$. Therefore, we have $L\le L'\le\alpha k$.

To see the fourth statement, note that there are at least $k^2/(2\alpha)$ fast
arcs by Claim~\ref{fast}, but fewer than $L'k/(4\alpha^2)\le k^2/(4\alpha)$
fast arcs start in points that are not good.

For the final statement, note that any $i$-fast arc has length at most $k+1$, by
definition. So, if such an arc starts at $q_j$, then it must be contained in
$c_j^*$. The curve $b_i$ must touch $c_j$ at the end point of the $i$-fast arc, because
the apex of the arc is $x_i$.
\qed
\medskip

We call an arc $z^*\subset b$ {\em short} if $|z^*\cap T|\le2\alpha k$. Note that while
the length counts sad touching points on an arc, in this definition we count
all touching points. We say that the
good points $q$ and $q'$ are {\em close}, if either the arc
of $b$ from $q$ to $q'$ or the arc from $q'$ to $q$ is short.

\begin{claim}\label{notclose}
Let $q$ and $q'$ be good points. If the arc $b^*$ from $q'$ to $q$ is short,
then $|b^*\cap X|\le2\alpha\alpha_1k$.

If $q_j$ and $q_{j'}$ are not close, then the arcs $c_j^*$ and $c_{j'}^*$ are
disjoint.
\end{claim}

\proof The first claim holds, because $q$ is a sad touching point, so
the arc $b^*$ ending there must have crossing-to-touching ratio below
$\alpha_1$.

For the second claim, assume that $c_j^*$ and $c_{j'}^*$ intersect and let $W$ be a
Jordan curve connecting $q_j$ to $q_{j'}$ along part of $c_j^*$ and
$c_{j'}^*$. Consider the two arcs $b^*$ and $b'^*$ that $b$ is cut  by $q_j$
and $q_{j'}$. By our assumption, neither of these arcs is short, so each has
more than $2\alpha k$ distinct curves of $\F$ touching it. As $W\cap X\le2\alpha k$, we
must have a curve $z\in \F$ touching $b^*$ that is disjoint from $W$.
Analogously, we have another curve $z'\in \F$ touching $b'^*$ and also
disjoint from $W$. Now $b$ and $W$ separate $z$
and $z'$, contradicting the fact that they (as any two curves in $\F$) must
intersect. This contradiction completes the proof of the claim. \qed

\medskip

We call the distinct good points $q_j$ and $q_{j'}$ {\em mingled} if
$|c_j^*\cap c_{j'}^*|>\alpha^2k/w$.

\begin{claim}\label{mingled}
Mingled points are close. A good point is mingled with at most $w/\alpha$
other good points.
\end{claim}

\proof The first statement follows directly from Claim~\ref{notclose}. The
second statement follows from the statement of Claim~\ref{goodpoor} that
$c_j^*$ contains at most $\alpha k$ intersection points in total. \qed
\medskip

For a good point $q$ let $I_q$ stand for the set of indices $1\le i\le k$ with
an $i$-fast arc starting at $q$. Similarly, for $1\le i\le k$, let $Q_i$ stand
for the set of good points $q$ at which an $i$-fast arc starts.

\begin{claim}\label{close} Let $1\le j<j'\le L$ be such that the arc from
$q_{j'}$ to $q_j$ is short, but $q_j$ and $q_{j'}$ are not mingled. Then
$|I_{q_j}\cap I_{q_{j'}}|<6\alpha^2k/\sqrt v$.
\end{claim}

\noindent{\it Proof sketch.} As the proof of Claim \ref{close} is fairly involved, we only sketch it here, while postponing the full details to Section \ref{Subsec:ProofClose}.

To simplify the presentation, let us first assume that the arc $c^*_j$ and $c^*_{j'}$ are disjoint. Denote $I=I_{q_j}\cap I_{q_{j'}}$.
Assume for a contradiction that $|I|\geq 6\alpha^2k/\sqrt v$.
The key observation is that at least ${|I|\choose 2}=\Omega(\alpha^4k^2/v)$ pairs of curves $b_{i},b_{i'}$ with $i,i'\in I$ determine  a lens of size $O(\alpha^3k)$ each.
As a result, $q_j$ receives at least $\alpha$ units of charge from such lenses, by the second rule, contrary to its choice as a poor point.

\begin{figure}[htbp]
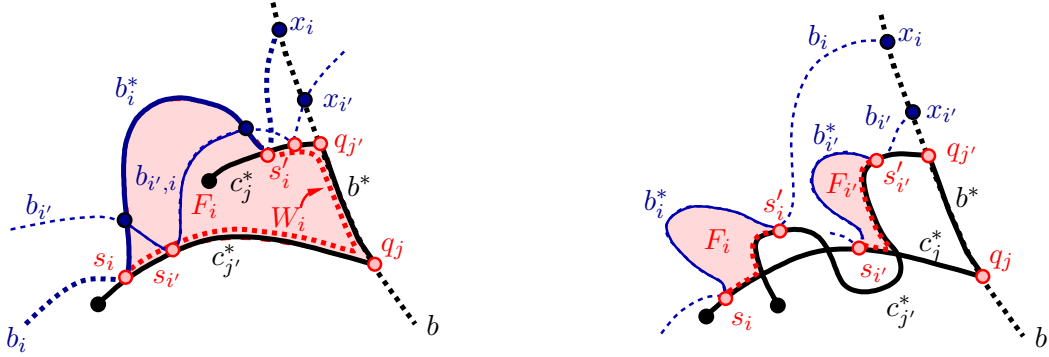

\begin{center}
\input{Lens1.pstex_t}\hspace{3cm}\input{Lens3.pstex_t}
\caption{\small Proof of Claim \ref{close}. Left: If $c_j^*$ and $c_{j'}^*$ are disjoint, then any pair of curves $b_i,b_{i'}$ with $i,i'\in I_{q_j}\cap I_{q_{j'}}$ yields a lens $b_{i',i}$ or $b_{i,i'}$. Right: If the arcs $c_j^*$ and $c_{j'}^*$ have multiple intersections, then a pair of curves $b_i$ and $b_{i'}$ can determine non-overlapping pockets $F_i$ and $F_{i'}$.  Thus, $b_i$ and $b_{i'}$ do not necessary determine a lens within $F_i$ or $F_{i'}$.}
\label{Fig:SketchClose}
\end{center}
\end{figure}

Indeed, let $b^*$ be the short arc of $b$ from $q_{j}$ to $q_{j'}$. Consider the curve $W=c^*_{j}\cup c'^*_j\cup b^*$ which contains at most $O(\alpha^2 k)$ intersection points with the curves of $\A$;
see Figure \ref{Fig:SketchClose} (left).
As each curve $b_i$, with $i\in I$, touches $W$ at a pair of points $s_i\in c^*_j$ and $s'_i\in c^*_{j'}$,  they determine a pocket $F_i$ which is bounded by (i) the portion $W_i$ of $W$ between $s_i$ to $s'_i$, and (ii) the arc $b_i^*$ of $b_i$ between the same two points and is to the right of $b_i^*$.
For any other curve $b_{i'}$ with $i'\in I\setminus\{i\}$, which touches $c_j$ within $W_i$, we define $b_{i',i}$ to be the shortest arc of $b_{i'}$ between two points of $b_i$ that contains $s_{i'}$. Note that $b_{i',i}$ is a lens inside $F_i$.

To bound the size of the lens $b_{i',i}$, we argue that each curve that touches $b_i^*$ must exit $F_i$ through $W_i$ (or, else, it will not meet one of the curves $c_j$ or $c_{j'}$). Hence, the overall number of such curves does not exceed $O(\alpha^2k)$. Since $b_i^*$ is adjacent to a sad point $s_i$, we obtain that $|b_i^*\cap X|=O(\alpha^3k)$.
Finally, each curve that touches the lens $b_{i',i}$ must also leave $F_i$ through $W_i\cup b_i^*$, so their number is also $O(\alpha^3k)$, in fact, at most $3\alpha^3k$. This means that, by the second charging rule, $b_{i',i}$ sends some charge to $q_j$.

The contradiction comes from $q_j$ being poor despite the fact that any pair of distinct indices $i,i'\in I$ determine a lens (either $b_{i',i}$ or $b_{i,i'}$) sending charge to $q_j$.

If $c_j^*$ and $c_{j'}^*$ intersect (possibly many times), then the above argument fails, as the curves $b_{i}$ may determine smaller-size pockets $F_i$ amidst $c_j^*\cup c_{j'}^*$, which do not necessarily overlap (see Figure \ref{Fig:SketchClose} (right)).
As a result, the number of lenses $b_{i',i}$ can be substantially smaller than ${|I|\choose 2}$, and it generally depends on the number of intersections between $c^*_j$ and $c^*_{j'}$.
Nevertheless, since $q_j$ and $q_{j'}$ are not mingled, the arcs $c^*_j$ and $c^*_{j'}$ have at most $\alpha^2k/w$ intersections, which enables to extend the previous analysis by finding somewhat fewer lenses sending charges to $q_j$ or $q_{j'}$, but noticing that these lenses tend to be shorter and therefore send more charge. For the precise accounting (see Section~\ref{Subsec:ProofClose}), we use Lemma~\ref{Lemma:RedBlue} again. \qed

\begin{claim}\label{fewclose}
Let $1\le j<j'\le L$ be such that the arc $b^*$ from $q_{j'}$ to $q_j$ is
short. The number of good points in $b^*$ is at most $50\alpha^3w$.
\end{claim}

\proof Let $S$ be the set of good points in $b^*$. We have
$\displaystyle |S|\frac k{4\alpha^2}\le\sum_{q\in S}|I_q|=\sum_{i=1}^k|Q_i\cap S|$
and, by the Cauchy--Schwarz inequality,
$\displaystyle \frac{k|S|^2}{16\alpha^4}\le\sum_{i=1}^k|Q_i\cap S|^2=\sum_{q,q'\in
S}|I_q\cap I_{q'}|.$
We use the trivial bound $|I_q\cap I_{q'}|\le k$ if $q=q'$ or if $q$ and $q'$ are
mingled. By Claim~\ref{mingled}, there are at most $(w/\alpha+1)|S|$ such
terms. The remaining terms can be bounded by $6\alpha^2k/\sqrt v$ using
Claim~\ref{close}. We obtain
$\displaystyle \frac{k|S|^2}{16\alpha^2}\le k(w/\alpha+1)|S|+\frac{6|S|^2\alpha^2k}{\sqrt v}.$
Substitute $v=21000\alpha^{12}$ and the claim follows. \qed.

\begin{figure}[htbp]
\begin{center}
\input{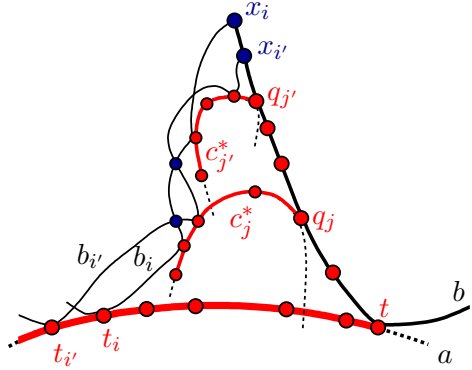}
\caption{\small Proof of Claim \ref{tower}.
The good points $q_j,q_{j'}\in S$ are not close, and their arcs $c^*_j$ and $c^*_{j'}$ form a pair of neighboring ``teeth''  in the comb $\Gamma$.}
\label{Fig:Tower}
\end{center}
\end{figure}

\begin{claim} \label{tower}
For distinct indices $1\le i,i'\le k$ the curves $b_i$ and $b_{i'}$ have at least
$|Q_i\cap Q_{i'}|/(100\alpha^3w)-1$ intersection points.
\end{claim}

\proof Let $Q=Q_i\cap Q_{i'}$. We select a subset $S\subseteq Q$ with no two
close points greedily: we consider the elements $q\in Q$ along $b$ in reverse
direction starting
from $t$, and include them in $S$ unless $q$ is close to a point already in
$S$. No short arc of $b$ avoiding $t$ contains  more than $50\alpha^3w$
good points by Claim~\ref{fewclose}. Thus, we have
$|S|\ge|Q|/(50\alpha^3w)-1$. Let $b^*$ be the
arc of $b$ from the point we put in $S$ last to the point $t$ we put there
first. By Claim~\ref{notclose}, the arcs $c_j^*$
corresponding to the points $q_j\in S$ are pairwise disjoint. Let $\Gamma$ be
the comb-like arrangement of these arcs together with $b^*$, see Figure \ref{Fig:Tower}. Let
$b_i^*$ be the maximal arc on $b_i$ from a touching point where an $i$-fast
arc ends to the apex $x_i$ of the $i$-fast arcs. Clearly, $b_i^*$ touches all
the ``teeth'' of the comb $\Gamma$, but it does not intersect its spine
$b^*$. This implies that $b_i^*$ touches the teeth in the same order as the
arc $b^*$. This is also true for the analogously defined arc $b_{i'}^*$ of
$b_{i'}$. Consider two neighboring teeth of the comb $\Gamma$. Obviously, either the
part of $b_i^*$ between the corresponding touching points is crossed by
$b_{i'}$, or the part of $b_{i'}^*$ between the corresponding touching points
is crossed by $b_i$. As the segments of $b_i^*$ between touchings of
consecutive teeth are disjoint, they represent at least $(|S|-1)/2$ intersections
between $b_i$ and $b_{i'}$.
\qed.

\medskip
Having proved these claims, we return to the proof of Lemma~\ref{in}. We started the proof by
assuming that the lemma fails, so we need to arrive at a contradiction to finish the proof.

By Claim~\ref{goodpoor}, we have many fast arcs starting at good points, namely
$ \sum_{j=1}^L|I_{q_j}|\ge\frac{k^2}{4\alpha}.$
Using $L\le\alpha k$ (Claim~\ref{goodpoor} again) and the Cauchy--Schwarz inequality, we get
$$ 
\frac{k^3}{16\alpha^3}\le\sum_{j=1}^L|I_{q_j}|^2=\sum_{1\le i,i'\le
k}|Q_i\cap Q_{i'}|.
$$
Subtracting the contribution of the case $i=i'$ and dividing by $2$, we obtain

\begin{equation}\label{Eq:Qii}
\sum_{1\le i<i'\le k}|Q_i\cap Q_{i'}|\ge\frac{k^3}{40\alpha^3}.
\end{equation}

Claim~\ref{tower} shows that this lower bound on the left-hand side of (\ref{Eq:Qii}) provides a lower bound on the
number of the intersection points between the curves $b_i$. We find that the number of such intersection points is at least
$k^3/(4000\alpha^6w)-k^2/2$. With the prior choice of parameters,
this contradicts the assumption that the total number of intersections
satisfies $|X|<\alpha_1n^2$. This contradiction proves Lemma~\ref{in}.

\subsection{Wrapping up the proof of Theorem \ref{main-bipartite}}\label{wrap} Finishing the proof of Theorem~\ref{main-bipartite} is simple once we have
Lemmas~\ref{happy}, \ref{out13}, \ref{out2}
and \ref{in}. Considering
all the charges in all the $M$ phases of our scheme, every sad touching point $t\in T'$ receives a charge of at
least $c_{\mbox{\scriptsize in}}=\alpha M$ by Lemma~\ref{in}. For an intersection
point $x\in X$, the total charge sent
out is at most $c_{\mbox{\scriptsize out}}=10M$ by Lemmas~\ref{out13} and \ref{out2}. Comparing the total charges sent and received we obtain
$$\frac{|X|}{|T'|}\ge\frac{c_{\mbox{\scriptsize in}}}{c_{\mbox{\scriptsize
out}}}=\alpha/10.$$
We have $|T'|\le 10|X|/\alpha$ from the line above and
$|T|-|T'|\le6|X|/\alpha_1$ from Lemma~\ref{happy}. In total, we have
$|T|\le16|X|/\alpha_1$, and the statement of the Theorem~\ref{main-bipartite} follows. \qed

\subsection{Proof of Claim
\ref{close}}\label{Subsec:ProofClose}

\proof For simplicity, we write $q$ and $q'$ for $q_j$ and $q_{j'}$,
respectively. Analogously, we write $c$, $c'$, $c^*$ and $c'^*$ for $c_j$,
$c_{j'}$, $c_j^*$ and $c_{j'}^*$, respectively. We write $b^*$ for the short arc
from $q'$ to $q$.

\begin{figure}[htbp]
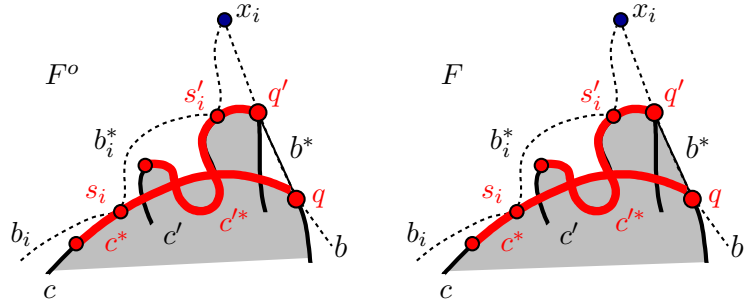

\begin{center}
\input{Faces.pstex_t}\hspace{1cm}\input{Faces1.pstex_t}
\caption{\small Proof of Claim \ref{close}. The curve $b$ lies in the single face $F^0$ of the arrangement of $c$ and $c'$.
The cell $F^0$ contains all the curves $b_i$ with $i\in I$ (left). The arc $b^*$ from
$q'$ to $q$ splits $F^0$ into two sub-faces, with all the touching points
$s_i$ and $s'_i$ lying on the boundary the sub-face $F\subset F^0$ to the left
of $b^*$ (right).}
\label{Fig:Faces}
\end{center}
\end{figure}

Refer to Figure \ref{Fig:Faces}.
Consider the arrangement of the curves $c$ and $c'$. The curve $b$
touches both of these curves, so it must be contained in a single face $F^0$ of
the arrangement, and this face is to the right of both $c$ and $c'$. As $c$ and
$c'$ intersect, the
boundary of $F^0$ is a simple closed Jordan curve which we denote by
$W^0$. Clearly, $W^0$ consists of alternating arcs of $c$ and $c'$ each
consistently oriented with $F^0$ to the right of them. The arcs of $c\cap W^0$ appear in the same cyclic order along $c$ and $W^0$, and a similar statement
is true for the segments of $c'\cap W^0$. Note, however, that outside $F^0$
the curves $c$ and $c'$ can behave wildly and all sorts of extra intersections
can occur even between $c^*$ and $c'^*$.

The arc $b^*$ splits this face in two, let $F$ stand for the
side of $F^0$ containing $b\setminus b^*$ (that is, on the left from
$b^*$). Let $W^1$ stand for the part of $W^0$ on the boundary of $F$. Clearly,
this is a Jordan curve connecting $q'$ to $q$ and all its segments coming
from $c$ and $c'$ are consistently oriented from $q'$ to $q$.

Let us write $I=I_q\cap I_{q'}$. In what follows, we can assume that $I$ is not empty.
For each $i\in I$, the curve $b_i$ touches both $c$ and $c'$ and intersects
$b$, so it is confined to the face $F^0$. Let $s_i\in c^*$ and
$s'_i\in c'^*$ be the points where $b_i$ touches the boundary of $F^0$. Recall
that $x_i$ is an intersection point of $b_i$ and $b$. Let $b_i^*$ be the arc
of $b_i$ from $s_i$ to $s'_i$ or vice versa,
whichever does not contain the apex $x_i$.

\begin{proposition}\label{Claim:ArcsBi}
With the previous notation the following holds.

\smallskip
(i) $b_i^*$ does not intersect $b$, and the apex
$x_i$ is the first point on $b_i$ after $b_i^*$ that belongs
to $b$.

(ii) We have $x_i\notin b^*$. (However, $b^*$ can still intersect $b_i$ at points other than $x_i$.)
\end{proposition}

\proof Part (i) follows from the definition of the apex $x_i$, and using that both $q$ and $q'$ are starting points of $i$-fast arcs.

For part (ii), recall that $q=q_j$ and $q'=q_{j'}$ are starting points of $i$-fast arcs, so both of them must lie in the interval of $b$ from $x_i$ to $t=q_1$.
As the good points $t=q_1,\ldots, q_j,\ldots, q_{j'},\ldots,q_L$ are listed in the reverse order of $b$, the four points $x_i,q_{j'},q_j,t=q_1$ must appear in this order along $b$. \qed

\medskip
Proposition \ref{Claim:ArcsBi} implies that the interior
of $b_i^*$ is inside the face $F$, and its endpoints $s_i$ and $s'_i$ are on the
boundary of $F$, and, consequently, also on $W^1$.
In what follows, we define a simple open Jordan curve $W$ that contains
$(c^*\cup c'^*)\cap W^1$. In order for $W$ to remain simple and connected, it may include segments outside $c^*$ and $c'^*$. Nevertheless, $W$ has at most $O(\alpha^2k)$ intersection points with the curves of $\A$.

\begin{figure}[htb]
\begin{center}
\input{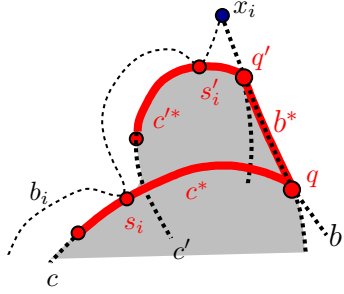}
\caption{\small Constructing the curve $W$. In the disjoint case, $W$ includes the arc $b^*$ from $q'$ to $q$.}
\label{Fig:TraceWDisjoint}
\end{center}
\end{figure}

\medskip
\noindent{\bf Tracing $W$.} For the definition of the curve $W$, we consider the following cases:

In case $c'^*$ does not intersect $c$, it lies entirely on $W^0$. In particular,  $c^*$ and $c'^*$ are disjoint. Hence, we can simply take $W$ to be
the union of $c^*$, $c'^*$ and $b^*$; see Figure \ref{Fig:TraceWDisjoint}. We call this the {\em the disjoint case}.

\begin{figure}[htb]
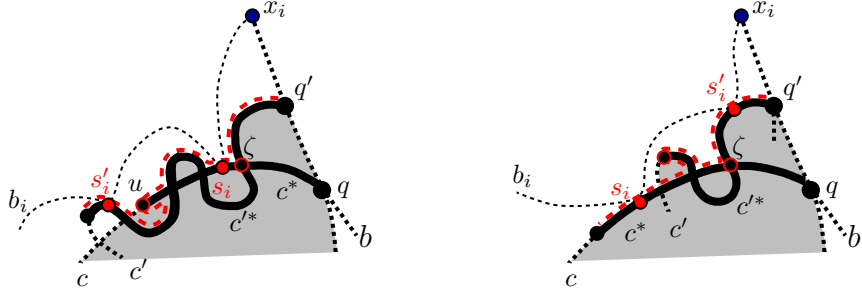

\begin{center}
\input{TraceIntersectingCase.pstex_t}\hspace{2cm}\input{TraceIntersectingFix.pstex_t}
\caption{\small Constructing the curve $W$ -- the two intersecting scenarios.
Notice that the first intersection $\zeta$ of $c'^*$ with $c$ must belong to $c^*$.  Left: In the first intersecting scenario, $W$ may have to retrace the last segment of $U$. Notice that $c'\cap W$ consists only of segments of $c'^*$ (while $c$ can meet $W$ at a point outside $c^*$). Right: In the second intersecting case, $W$ may have to retrace the last segment of $U'$, and $c\cap W$ consists only of segments of $c^*$.}
\label{Fig:TraceWIntersect}
\end{center}
\end{figure}

In case $c'^*$ intersects $c$, we consider the first such intersection point $\zeta$ along $c'^*$, at which $c'^*$ leaves $W^0$, and notice that $\zeta$ must belong to $c^*$; see Figure \ref{Fig:TraceWIntersect}. Indeed, assume for a contradiction that $\zeta$ belongs to $c\setminus c^*$.
Since (i), the order of the segments of $c\cap W^0$ along $W^0$ is consistent with their order along $c$, (ii) $W^1$ ends at the starting point $q$ of $c^*$, and (iii) $W^1$ begins at $q'$,
the last appearance of $c^*$ along $W^0$ is also contained in $W^0\setminus W^1$. However, in that case, $c^*$ can never show
up on $W^1$, contrary to $I\neq \emptyset$ (and, thus, to the choice of $s_i$ on $W^1\cap c^*$ for $i\in I$). We can assume, then that the first intersection $\zeta$ of $c'^*$ with $c$ lies on $c^*$, and it is the first appearance of $c$ and $c^*$ along $W^1$.

As $c'^*$
starts at $q'$, there is a shortest segment $U'$ of $W^1$ starting at $q'$ with $c'^*\cap
W^1=c'\cap U'$. Although $c^*$ starts at $q$, we similarly define the shortest interval $U\subset W^1$ starting at $q'$ with the property that $c^*\cap
W^1=c\cap U$. Indeed, $c^*\cap W^0$ is
confined to an interval of $W^0$ starting at $q$ and we can choose $U$ to be
the intersection\footnote{This intersection is indeed a connected interval of $W^1$ since $W^1$ ends at $q$.} of the shortest such interval with $W^1$.

We call the case when $U$ is shorter than $U'$ the {\em first intersecting case} (see Figure \ref{Fig:TraceWIntersect} (left)). In this case, we have $c^*\cap
W^1\subseteq U\subseteq c^*\cup c'^*$. If the end point $u$ of $U$ is in
$c'^*$, we let the curve $W$ start with $U$ followed by the part of $c'^*$ after $u$. If
$u\notin c'^*$, then we let $W$ start with $U$ again, but we cannot directly add a part of $c'^*$, so instead we follow $U$ by a curve retracing the last
segment of $U$ (which is contained in $c^*$) very close, but slightly outside $F$ till we meet $c'^*$, and
then add the remaining part of $c'^*$.
Assuming that the last segment of $U$ is retraced sufficiently close to $c^*$, any curve
that meets the ``reversed'' segment $\gamma$ must also meet $c^*$.

In case $U'$ is shorter than $U$ (the {\em second intersecting case}),
we construct $W$ symmetrically: We start with $U'$ and add the remaining part of $c^*$,
possibly using a reverse segment $\gamma$ slightly outside $F$ to connect the two
parts. See Figure \ref{Fig:TraceWIntersect} (right).

Note that the only case when this construction does not work is when $U'$ consists of a single segment so we could follow it backward all the way and still not meet $c$. But this cannot happen, because configurations like that are treated separately in the disjoint case.

The properties of $W$ are summarized in the following proposition.

\begin{proposition}\label{Claim:TraceW}
$W$ is a simple open Jordan curve. It contains $(c^*\cup c'^*)\cap W^1$ and consists only of segments of $c^*$, $c'^*$, and possibly one additional segment.
Specifically, $W$ satisfies the following properties:
\begin{itemize}
\item[(a)] In the disjoint case (when $c'^*$ does not intersect $c$), we have $W=c^*\cup c'^*\cup b^*$, $b\cap W=b^*$ and $c\cap W=c^*$.
\item[(b)] In both intersecting cases,
(b1) $W$ lies outside the interior of $F$, and is composed of segments of $c^*$ and $c'^*$, and possibly of an additional
segment $\gamma$ closely retracing\footnote{Namely every intersection of $\gamma$ with a curve $\sigma \in \A$ corresponds to a unique crossing (non-touching intersection) of $\sigma$ with the retraced segment of $W$ and vice versa.
$\gamma$ does not not contain points of $X\cup T$.} the previous segment of $W$ from the outside of $F$, (b2) the order of the segments of $c^*\cap W$ and $c'^*\cap W$ along $W$ is consistent with the respective orientations of $c$ and $c'$, and (b3) $b\cap W=\{q'\}$.

\item[(c)] In the first intersecting case, $c'\cap W$ consists of segments of $c'^*$. When following $W$ past the end point of any of these segments, the curve $W$ continues on the right of $c'$.

\item[(d)] In the second intersecting case, $c\cap W$ consists of segments of $c^*$. When following $W$ past the end point of any of these segments, the curve $W$ continues on the right of $c$.
\end{itemize}
\end{proposition}

\proof
Most of the statements follow directly from our construction of $W$. Note, however, that in the intersecting case $c$ and $c'$ or even $c^*$ and $c'^*$ can intersect in unexpected ways outside (the closure of) $F^0$. The curve $W$ is still simple, as it consists of a segment (namely, either $U$ or $U'$) of the (simple) boundary of $F$ followed by a segment of the simple curve $c^*$ or $c'^*$, possibly with a ``retracing'' curve $\gamma$ in between. In the first intersecting case, $c'^*$ contains this final segment, so the remaining part of $c'\cap W$ is all from $c'^*\cap W^1$. Property (b3) follows from (b2) and since $c^*$ is a proper subarc of $c$, so none of the segments of $c^*\cap W^1$ can be adjacent to $q$. Part (c) holds as $F^0$ is to the right of all the segments on its boundary. Note, however, that $c\cap W$ may contain several ``unintended'' intersection points of $c$ and (this last part of) $c'^*$. Furthermore, $\gamma$ (if it exists) is to the left of the segment it retraces. Part (d) can be seen similarly, with the roles of $c$ and $c'$ reversed. \qed

\medskip
Let us define $B$ as the set of points where $W$ is intersected by a
curve in $\A$ (other than the curve $W$ follows at that segment). We have $X\cap W\subseteq B$, but $B$ may contain further
crossing points along the reverse segment $\gamma$, if such a segment exists. Still, in the intersecting case we have $|B|\le3\alpha k$ by Claim~\ref{goodpoor}(iii), while $|B|\le2\alpha\alpha_1k+2\alpha k$, where we use also Claim~\ref{notclose} to estimate the size of $B\cap b^*=X\cap b^*$. We will use the bound $|B|\le2\alpha^2k$ that holds in both cases and
apply Lemma~\ref{Lemma:RedBlue} to $B$ with a non-uniform weight function
$w_0$. We set $w_0(x)=w+1$ for $x\in B\cap c^*\cap c'^*$ or (in the disjoint
case) if $x=q$. We set $w_0(x)=1$ otherwise. For the total weight we have
$w_0(B)\le3\alpha^2k$ as $q$ and $q'$ are not mingled.

For $i\in I$, let $W_i$ be the portion of $W$ between the touching points $s_i$ and $s'_i$,
and let us write $l_i=w_0(B\cap W_i)$. Let $R_0=\{s_i\mid i\in I\}$ and set
$\lambda=3\sqrt v$. Let
$R=\{s_i\mid i\in I, |W_i\cap R_0|\le l_i/\lambda\}$ and $\hat I=\{i\in I\mid
s_i\notin R\}$. We have $|I|=|R|+|\hat I|$. In what follows, we bound $|R|$ and
$|\hat I|$ separately.

\medskip
\noindent{\bf Bounding $|R|$.} We apply Lemma~\ref{Lemma:RedBlue} for the curve $W$, the sets
$R$ and $B$, the weight function $w_0$ and the parameter $\lambda$. The
condition is satisfied, as for the interval $W_i$ ending at
$s_i\in R$ we have $|W_i\cap R|\le|W_i\cap R_0|\le l_i/\lambda=w_0(W_i\cap
B)/\lambda$. From Lemma \ref{Lemma:RedBlue}, we conclude that
\begin{equation}\label{Rbound}
|R|\le\frac{3w_0(B)}\lambda\le\frac{9\alpha^2k}\lambda.
\end{equation}

\begin{figure}[htbp]
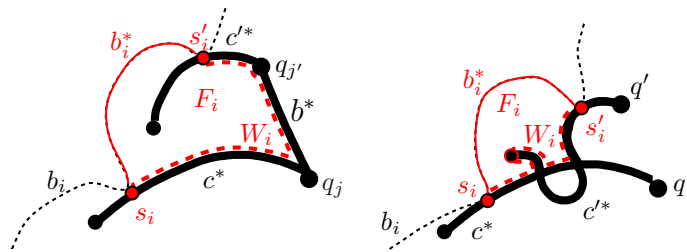

\begin{center}
\input{TraceDisjoint2Fix1.pstex_t}\hspace{0.8cm}\input{TraceIntersecting2Fix.pstex_t}
\caption{\small The portion $W_i$ of $W$ between $s'_i=b_{i}\cap c'$ and $s_i=b_i\cap c$ is traced.
The face $F_i$ of $\reals^2\setminus (b_i\cup W_i)$ lies to the right of $b_i$. Left: The arcs $c^*$ and $c^{\prime*}$ are disjoint so both $W$ and $W_i$ must include $b^*$.  Right: The scenario where $c^*$ and $c^{\prime*}$ intersect.}
\label{Fig:TraceW1}
\end{center}
\end{figure}

\medskip
\noindent{\bf Bounding $|\hat{I}|$.} Fix $i\in \hat{I}$, whose respective arc $W_i$ contains at least $l_i/\lambda$ points $s_i\in R_0$.
Consider $W_i\cup b_i^*$ and let $F_i$ be the
side of this closed Jordan curve to the right of $b_i^*$; see Figure \ref{Fig:TraceW1}.

\begin{figure}[htb]
\begin{center}
\input{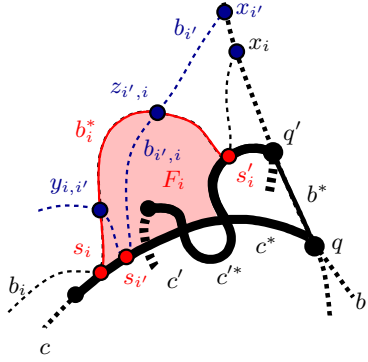}
\caption{\small The segment $W_i$ contains a touching $s_{i'}$ between $c$ with $b_{i'}$, for $i'\in I\setminus \{i\}$. Each such $b_{i'}$ yields a lens $b_{i',i}$ of length at most $\alpha(l_i-w)$.}
\label{Fig:TraceLens}
\end{center}
\end{figure}

First we sketch our argument for our bound on $|\hat{I}|$. The rigorous calculation is after Propositions~\ref{ClaimBiShort} and \ref{Claim:DisjointFi}. We will show that any curve $b_{i'}$ whose respective touching $s_{i'}$ lies on $W_i$, determines within $F_i$ a lens $b_{i',i}$ of length at most $3\alpha^3k$; see Figure \ref{Fig:TraceLens}.
Each of these lenses will send some charge to $q$ via the second rule. The final bound on $|\hat I|$ follows from $q$ being poor despite all these incoming charge.

The crux of our argument is showing that each lens $b_{i',i}$ has length at most $\alpha(l_i-w)$. To this end, we first bound the quantities $|b_i^*\cap T|$ and $|b_i^*\cap X|$ in terms of the overall weight $l_i$ of $W_i$.

\begin{proposition}\label{ClaimBiShort}
With the previous assumptions, $b_i^*$ contains at most
$l_i-w$ points of $T$, and at most $\alpha_1(l_i-w)$ points of $X$.
\end{proposition}

The proof of Proposition \ref{ClaimBiShort} relies on the following topological property. We will use this observation to argue that no curve in $\cal A$ can stay inside $F_i$. Indeed, curves in $\cal G$ have to intersect $b$ and curves in $\cal F$ have to intersect both $c$ and $c'$.

\begin{proposition}\label{Claim:DisjointFi}
The curve $b$ does not intersect the interior of $F_i$. Furthermore, at least
one of the curves $c$ and $c'$ does not intersect the interior of $F_i$.
\end{proposition}

\proof
The proof is based on controlling where $b$ and one of $c$ and $c'$ intersects the boundary $W_i\cup b_i^*$ of $F_i$.

In the disjoint case, $b^*$ separates along $W$ the points $s_i$ on $c^*$ and $s'_i$ on $c'^*$, so  $b^*\subset W_i$. This means $b\cap W_i=b^*$ by Proposition~\ref{Claim:TraceW}(a). Note that the curve $b$ continues outside $F_i$ after it leaves the boundary of $F_i$ at $q$. Since $b\cap b_i^*=\emptyset$ by Proposition~\ref{Claim:ArcsBi} (and $b$ cannot cross $c^*\cap W_i$ or $c'^*\cap W_i$), the arc $b\setminus b^*$ too cannot enter the interior of $F_i$. We similarly argue that $c\setminus c^*$ cannot enter the interior of $F_i$ through either of the boundary arcs $b_i^*$ and $c'^*\cap W_i$.

For the rest of the proof, we consider the intersecting cases. We still have $b\cap b_i^*=\emptyset$ by Proposition~\ref{Claim:ArcsBi}. Now we have $b\cap W=\{q'\}\notin W_i$ by Proposition~\ref{Claim:TraceW}(b), so $b$ never intersects the boundary of $F_i$. Following $W$ past $s_i$ and $s'_i$, we see that it leaves $W_i$ and $F_i$ in one of these directions and reaches the point $q'\in b\setminus F_i$. This implies that $b$ is again outside $F_i$.

In the first intersecting case, we consider the curve $c'$. It touches $b_i$, so we have $c'\cap(W_i\cup b_i^*)=c'\cap W_i\subseteq c'\cap W$ and this is covered by Proposition~\ref{Claim:TraceW}(c). Recall that the part of $c'$ on the boundary of $F_i$ consists of segments of $c'^*$. It is easy to see that $F_i$ is on the right of each of these segments, just as $W_i$ continues to the right of $c'$ after the end points of any of these segments (again, by Proposition~\ref{Claim:TraceW}(c)). This shows that $c'$ never enters the interior of $F_i$.

A similar argument in the second intersecting case shows that $c$ never enters the interior of $F_i$. \qed

\medskip
\noindent{\bf Proof of Proposition \ref{ClaimBiShort}.}
The bound on $|b_i^*\cap T|$ follows from the fact that any curve $\sigma\in \F$
touching $b_i^*$ must intersect $W_i$.
Indeed, any such curve $\sigma$ is in $F_i$ in a
small neighborhood around the point where it touches $b_i^*$.
Since the curve $\sigma$ intersects {\it each} of the curves $c,c'\in \F$ (and at least one of $c$ and $c'$ is disjoint from the interior of $F_i$), $\sigma$ must meet $\partial F_i$ at a point of $B\cap W_i$.
We have $|B\cap W_i|\le l_i-w$, as $W_i$ contains at
least one of the heavy points with weight $w+1$. The bound on $|b_i^*\cap X|$
follows since the endpoint $s_i$ of $b_i^*$ is sad. \qed

\medskip
We are now ready to bound the cardinality of $\hat{I}$ and thus complete the proof of Claim \ref{close}.
Let us consider $i,i'\in I$ with $i\ne i'$ and $s_{i'}\in W_i$.
Follow $b_{i'}$ from $s_{i'}$ in both directions. It starts out inside $F_i$
and eventually has to reach $b$ that is disjoint of the interior of $F_i$. As
the part of $b_{i'}$ around $s_{i'}$ is in $F$, the first intersection with
$b$ in either direction is outside $b^*$.  Proposition \ref{Claim:DisjointFi} implies that $b_{i'}$ must properly cross the boundary of $F_i$ to meet $b$. The first intersection point in
either direction must be on $b_i^*$, for it cannot be on $b^*$, and $b_{i'}$
touches both $c_j$ and $c_{j'}$. Let us call these intersection points $y_{i',i}$
and $z_{i',i}$ such that the arc $b_{i',i}^*$ from $y_{i',i}$ to $z_{i',i}$ along
$b_{i'}$ is a lens, is inside $F_i$ and it contains $s_{i'}$; see Figure \ref{Fig:TraceLens}.

Note that $b_{i',i}^*$ is a lens. The length $l_{i',i}$ of this lens is at most
$\alpha(l_i-w)$. Indeed, Proposition \ref{Claim:DisjointFi} implies that any curve touching $b_{i',i}^*$
must intersect the boundary of $F_i$ so as to be able to intersect {\it both} $c$ and $c'$. However, by Proposition \ref{ClaimBiShort}, the
total number of points on this boundary at which a curve of $\A$ intersects it, is
at most $|B\cap W_i|+|b_i^*\cap X|\le(\alpha_1+2)(l_i-w)=\alpha(l_i-w)$.

Each such lens $b_{i',i}$ sends a
charge of $v/((l_{i',i}+w)k)$ to $q$, by the second
charging rule. Indeed, this rule applies to $b_{i',i}$, because its length is
$l_{i',i}\le\alpha l_i\le\alpha|B|\le 3\alpha^3k$, and
the arc from $q$ to $s_{i'}$ satisfies the requirements. The amount of the
charge sent is
$$\frac v{(l_{i',i}+w)k}\ge\frac v{(\alpha(l_i-w)+w)k}\ge\frac
v{\alpha(l_i-\lambda)k}.$$

If we further assume that $i\in\hat I$, then we have more than $l_i/\lambda$
choices of $i'\in I$ with $s_{i'}\in W_i$. One of these choices is $i'=i$, but
more than $l_i/\lambda-1$ other choices will give rise to lenses $b_{i',i}$,
each sending a charge of at least $v/(\alpha(l_i-\lambda)k)$ to
$q$. The total of these charges for a fixed $i\in\hat I$ is at least
$v/(\alpha\lambda k)$, and for all $i\in\hat I$ this is at least
$|\hat I|v/(\alpha\lambda k)$.

We know that $q_j$ is poor, so this charge does not reach the threshold of
$\alpha$. As a consequence, we have
$\displaystyle |\hat I|\le\alpha^2\lambda k/v.$
To finish the proof of Claim~\ref{close}, we use this last estimate, Equation~(\ref{Rbound}),
the fact $|I|=|R|+|\hat I|$ and substitute $\lambda=3\sqrt v$. \qed

\section{Proof of Theorem \ref{general}}\label{Sec:Dense}

To prove Theorem~\ref{general}, we have to get rid of the assumption in Theorem~\ref{main} that the curves are pairwise intersecting. We achieve this in two steps. First, in Subection~\ref{firsthalf}, we state and prove a separation result between the number of touchings and the number of intersections that does not assume strict pairwise intersection, but still assumes a very dense intersection graph. Then, in Subection~\ref{secondhalf}, we apply planar separation arguments to get rid of this milder assumption. In both of these steps, we lose in the separating function, namely, the exponent of the $\log\log$ function decreases. We did not attempt to optimize for this exponent, because we believe that even a logarithmic separation should hold, as stated in Conjecture~\ref{conjecture}.

\subsection{Sampling}\label{firsthalf}

In this subsection, we prove the following lemma. Note that, like Theorem~\ref{general}, it is about open Jordan curves, not closed ones.

\begin{lemma}\label{denselemma}
Let $\A$ be a family of $n$ simple open Jordan curves in general position in the plane. Let $T$ be the set of touching points between curves of $\A$ and let $X$ be the set of intersection points. With $h=n^2/|T|$ and $f=|X|/|T|$ we have $f^{72}h^{144}=\Omega(\log\log n)$.
\end{lemma}

\proof
As the statement of the lemma is asymptotic, we may assume in the calculations below that $n$ is sufficiently large.

\begin{figure}[htbp]
\begin{center}
\input{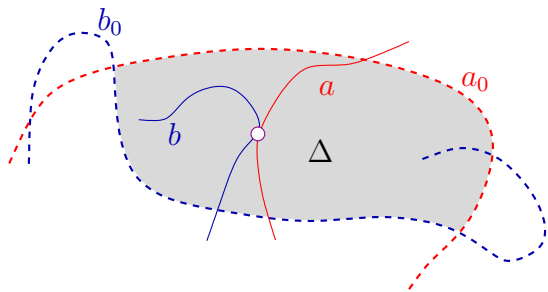}
\caption{\small Proof of Lemma~\ref{denselemma}. We select a pair $a_0,b_0\in \A$ satisfying three properties including that at least $|T|/(20h^2)$ of the touchings involve a pair of curves $a,b\in \A\setminus\{a_0,b_0\}$ both intersecting $a_0\cup b_0$. We then choose an open cell $\Delta\subset\reals^2\setminus(a_0\cup b_0)$ which contains at least $|T|/(2400fh^3)$ of these touchings.}
\label{Fig:Sampling}
\end{center}
\end{figure}

We select a pair of distinct curves $a_0,b_0\in \A$.  We try to select them so as to satisfy these conditions (see Figure \ref{Fig:Sampling}):

(a) For $m_0=|a_0\cap b_0|$ we want $m_0\le120fh$.

(b) For $m_1=|(a_0\cup b_0)\cap X|$ we want $m_1\le80nfh$.

(c) Let $m_2$ be the number of touchings $t\in T$ between two curves in $\A\setminus\{a_0,b_0\}$ with both of these curves intersecting $a_0\cup b_0$ (see Figure \ref{Fig:Sampling}). We want $m_2\ge|T|/(20h^2)$.

By selecting the pair $(a_0,b_0)$ uniformly at random, the expectation of $m_0$ is $E[m_0]=|X|/{n\choose2}<3f/h$. We have $E[m_1]\le2|X|/n=2nf/h$. For the expectation of $m_2$ notice that at most $|T|/2$ touchings are contained in a curve $a\in\A$ with $|a\cap T|\le|T|/(2n)$. The remaining elements of $T$ (at least $|T|/2$ of them) are counted in $m_2$ with probability at least $|T^2|/(5n^4)$ each, yielding $E[m_2]\ge|T|^3/(10n^4)=|T|/(10h^2)$.

By Markov's inequality condition (a) fails with probability less than $1/(40h^2)$ and the same holds for condition (b). Using Markov's inequality again and the fact that $m_2\le|T|$ we have that condition (c) is satisfied with probability at least $1/(20h^2)$. Thus, all three conditions are simultaneously satisfied with some positive probability. We fix such a choice of the curves $a_0$ and $b_0$ and call them the {\em ground curves}.

Let the family $\A'\subseteq\A\setminus\{a_0,b_0\}$ consist of the curves that intersect at least one of $a_0$ or $b_0$. By property (c), these curves create $m_2\ge|T|/(20 h^2)$ touchings.

 Let us consider the arrangement of $a_0$ and $b_0$. In case the ground curves are disjoint there is a single cell of this arrangement and its boundary is $a_0\cup b_0$. We will treat this somewhat peculiar case later. Otherwise, the arrangement has $m_0\le120fh$ cells, each with a connected boundary. For each open cell $\Delta\subset \reals^2\setminus (a_0\cup b_0)$, let $T_{\Delta}\subset T\cap \Delta$ denote the set of touching points between the curves of $\A'$ within $\Delta$ (again, see Figure \ref{Fig:Sampling}). By the pigeon-hole principle,
there exists an open cell $\Delta\subset \reals^2\setminus (a_0\cup b_0)$ with
\begin{equation*}
|T_\Delta|\geq |T|/(2400fh^3).
\end{equation*}

\begin{figure}[htbp]
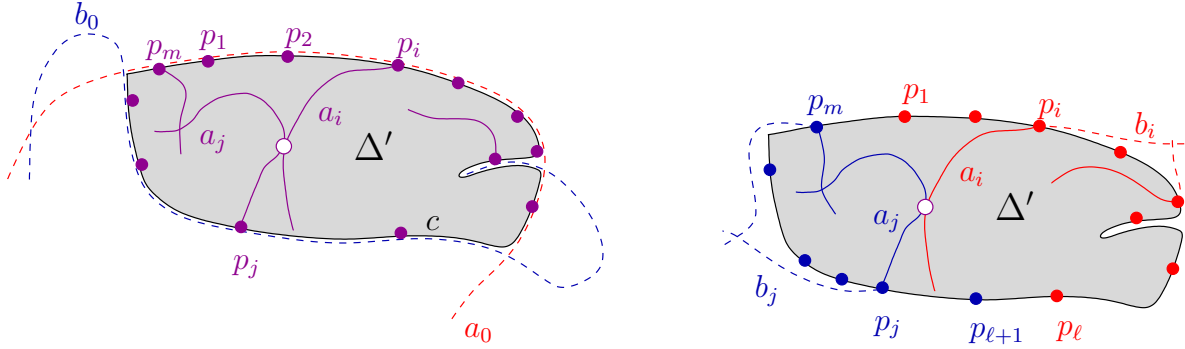

\begin{center}
\input{SamplingFix1.pstex_t}\hspace{1.5cm}\input{SamplingFix2.pstex_t}
\caption{\small Proof of Lemma~\ref{denselemma} -- constructing the families $\F$ and $\G$. Left: We obtain a new family $\A''$ by trimming each curve of $\A$ to $\Delta$. We then slightly inflate the boundary of $\Delta$ to a simple closed Jordan curve $c$ which encloses a region $\Delta'\subset \Delta$ and meets each $a_i\in \A''$ at a single point $p_i$. Right: We choose a random index $1\leq \ell \leq m$ and augment each $a_i\in \A''$ with a curve $b_i$ outside $\Delta'$, with the property that $b_i$ and $b_j$ intersect at a single point if and only if $1\leq i,j\leq \ell$ or $\ell<i,j\leq m$, and otherwise they are disjoint.}
\label{Fig:Sampling1}
\end{center}
\end{figure}

We fix such a cell $\Delta$. We consider each connected component of $a\cap\Delta$ for curves $a\in\A'$. These are simple Jordan curves in $\Delta$ with at least one end point on the boundary. We make the curves slightly shorter to make sure each has exactly one endpoint on the boundary but they still determine the same set $T_\Delta$ of touchings. We denote by $\A''$ the resulting family of $m\le n+m_1\le(80fh+1)n$ curves; see Figure \ref{Fig:Sampling1} (left).

We can slightly inflate the boundary of $\Delta$ to a simple closed Jordan curve $c\subset\Delta$ with $c$ intersecting each curve $a\in\A''$ exactly once, close to the end point of $a$ on the boundary of $\Delta$ and with all the touching points $T_\Delta$ on the side $\Delta'$ of $c$ contained in $\Delta$. Let us enumerate the curves in $\A''$ as $\A''=\{a_1,a_2,\dots,a_m\}$ such that the intersection points $p_i=a_i\cap c$ appear on $c$ in this cyclic order.

We pick a parameter $1\le\ell\le m=|\A''|$ and form the family $\F$ by slightly modifying the curves $a_i$ for $1\le i\le\ell$ and form $\G$ by slightly modifying $a_i$ for $\ell<i\le m$. The slight modification consists of keeping $a_i\cap\Delta'$ and attaching a curve $b_i$ to it that starts at $p_i$ and is disjoint from $\Delta'$; see Figure \ref{Fig:Sampling1} (right). We choose these additional curves such that (i) the curves $b_i$ and $b_{i'}$ are disjoint if $i\le\ell<i'$, (ii) the distinct curves $b_i$ and $b_{i'}$ intersect exactly once if $i,i'\le\ell$ or $i,i'>\ell$ and (iii) the curves $b_i$ are in general position. Clearly, such curves $b_i$ exist.

The family $\F$ consist of at most $m$ pairwise intersecting Jordan curves and the same is true for $\G$. With adding dummy curves we can actually assume that both families consist of $m+n$ curves.

Let us choose $\ell$ uniformly at random. The touching between the curves $a_i$ and $a_{i'}$ with $i<i'$ will remain a touching between the corresponding modified curves (one in $\F$, the other one in $\G$) if we have $i\le\ell<i'$. This happens with probability $(i'-i)/m$. Among the $|T_\Delta|\ge|T|/(2400fh^3)$ such touchings at most half can be between curves $a_i$ and $a_{i'}$ with $i<i'<i+x$ for $x=|T|/(4800fh^3m)$. Each touching point of the other half of $T_\Delta$ remains a touching points between a curve in $\F$ and a curve in $\G$ with probability at least $x/m$. Thus, the expected number of touchings between a curve in $\F$ and  curve in $\G$ is at least $(|T_\Delta|/2)\cdot(x/m)$. We choose $\ell$ such that the actual number of these touchings is at least this expectation.

Some of the intersection points between curves in $\F\cup\G$ come from $X$. We have one remaining intersection point between any two curves in $\F$ and also between any two curves in $\G$ for a total of $O(m^2+n^2)$ additional intersection points.

We are almost ready to apply Theorem~\ref{main-bipartite} to finish the proof. The only hurdle yet to clear is to pass from open Jordan curves to closed ones. This can simply be done by slightly inflating each curve. The process can be done in such a way that (i) touching curves remain touching, (ii) intersecting curves remain intersecting, (iii) the general position property is preserved and (iv) the number of intersection points is multiplied by at most $4$. With this we obtain families $\F'$ and $\G'$, each consisting of $m+n$ pairwise intersecting simple closed Jordan curves with the total number of intersections between curves in $\F\cup\G$ being $O(|X|+m^2+n^2)=O(f^2h^2n^2)$ and with the number of tangencies between a curve in $\F$ and a curve in $\G$  at least $(|T_\Delta|/2)\cdot(x/m)=\Omega(f^{-4}h^{-10}n^2)$. Applying Theorem~\ref{main-bipartite} to these families $\F'$ and $\G'$ yields the statement of the lemma.

Finally, we have to consider the special case when the ground curves $a_0$ and $b_0$ chosen in the first step of our proof are disjoint. In this case the arrangement of the ground curves has a single cell $\Delta=\reals^2\setminus (a_0\cup b_0)$. We define $\A''$ exactly as in the general case, so $\A''$ consists of curves contained in $\Delta$ with one end point on one of the ground curves. We distinguish ``type-a'' or ``type-b'' curves in $\A''$ according to whether it has an end point on $a_0$ or on $b_0$.

If at least one third of the touchings between two curves of $\A''$ are between two type-a curves, then we simply ignore $b_0$ and the type-b curves and consider the cell $\Delta^*=\reals^2\setminus a_0$ and the type-a curves. We can finish the proof as in the general case as $\Delta^*$ has a connected boundary.

If at least one third of the touchings between two curves of $\A''$ are between two type-b curves, then we proceed analogously.

If none of the above two cases hold, then we concentrate on the touchings between a type-a and a type-b curve: at least one third of all touchings between curves of $\A''$ must be like this. The situation is even simpler in this case with no need for any random choice. We obtain $\F$ by modifying slightly  the type-a curves and we obtain $\G$ by modifying slightly the type-b curves. For this we have to separately inflate the two ground curves. We finish the proof as in the general case. \qed

\subsection{Separation}\label{secondhalf}

The main ingredient we need for the proof of Theorem~\ref{general}
is the following separator theorem of Fox and Pach \cite{FP08} for {\it intersection graphs} of families of Jordan arcs.

\begin{theorem}\label{Thm:Separator}
For any collection $\A$ of $n$ Jordan arcs in the plane in general position with a total of $x$ intersection points, there is a subset $\B\subseteq\A$ of cardinality $O(\sqrt{x})$ so that $\A\setminus \B$ can be divided into disjoint subsets $\A_1$ and $\A_2$ of cardinality at most $2n/3$ each with the property that no arc of $\A_1$ meets an arc of $\A_2$.
\end{theorem}
\medskip

Repeated application of this result yields the following.

\begin{lemma}\label{densepart}
Let $\A$ be a collection of $n$ simple Jordan curves in general position in the plane. Let $d$ be maximum number of intersection points in a curve of $\A$ and let $T$ be an arbitrary subset of the intersection points with $|T|^2\ge nd^3$. There exist a subset $\A_0\subseteq\A$ of $\Theta(n^2d^3/|T|^2)$ curves such that $\Omega(nd^3/|T|)$ of the points in $T$ are intersection points of curves of $\A_0$.
\end{lemma}

\proof Let us set a threshold parameter $1\le M\le n$. We split the $\A$ by finding a small separator subset $\B\subset \A$ and partitioning $\A\setminus\B$ into $\A_1$ and $\A_2$,  as described in Theorem~\ref{Thm:Separator}. We recursively apply this procedure to the families $\A_1$ and $\A_2$, stopping only when we obtain subsets of size less than $M$. Let $\B'$ denote the set of all separator arcs that are removed at any invocation of the recursion and let $\A'_i$ for $i\in I$ stand for the final partition of $\A\setminus\B'$.  Note that, by the properties of the separation, the curves in $\A'_i$ do not intersect curves from other parts $\A'_j$, $j\ne i$.

When we split a set $\A'\subseteq\A$ of size $|\A'|=m\ge M$, we find a subset $\B'\subseteq\A'$ with size $|\B'|=O(\sqrt{dm})$, as the curves in $\A'$ determine at most $dm$ intersection points. The resulting parts of $\A'\setminus\B'$ have size at most $2m/3$.

The final parts produced by the above partitioning satisfy $|\A'_i|<M$, since this was our halting condition. Taking into account that these parts were obtained by splitting a subset of size at least $M$, they also satisfy $|\A'_i|\ge M/3-O(\sqrt{Md})$.

Consider all the different parts $\A'$ obtained in intermediate steps of the above partitioning process. Those sets of size $|\A'|=m$ in an interval $(3/2)^iM\le m<(3/2)^{i+1}M$ for some fixed integer $i$ are clearly pairwise disjoint, so their number is at most $n/((3/2)^iM)$. We find a separator of size $O(\sqrt{(3/2)^iMd})$ for each one of them. The total contribution of this interval to the size of $\B'$ is $O((3/2)^{-i/2}M^{-1/2}nd^{1/2})$. We sum these contributions for the integers $0\le i\le\log_{3/2}n$, and obtain $|\B'|=O(n\sqrt{d/M})$.

We set $M=Cn^2d^3/|T|^2$ with a large constant $C$. We may assume this yields $M\le n$, for otherwise the choice $\A_0=\A$ satisfies the requirements of the lemma. It is clear from our bound on $|\B'|$ that if $C$ is large enough, we have $|\B'|\le|T|/(2d)$. Analogously, from $|\A'_i|\ge M/3-O(\sqrt{Md})$ that holds for all $i\in I$, we get $|\A'_i|\ge M/4$ if $C$ is large enough.

During the partition, we lose at most $|\B'|d\le|T|/2$ intersections between
the curves of $\A$, that is, for at least half of the points $t\in T$, neither
of the two curves of $\A$ through $t$ are in $\B'$. Both of these curves are
therefore in the same part $\A'_i$ as curves from distinct parts are disjoint. By the pigeonhole principle, there is a part $\A'_{i_0}$ with at least $|T|/(2|I|)$ of the points of $T$, showing up as intersection points between curves of $\A'_{i_0}$. The choice $\A_0=\A'_{i_0}$ satisfies the requirements of the lemma since $|A'_{i_0}|=\Theta(M)=\Theta(n^2d^3/|T|^2)$ and $|T|/(2|I|)=\Omega(nd^3/|T|)$, where we use that $|A'_i|\ge M/4$ and therefore $|I|\le n/(M/4)=O(|T|^2/(nd^3))$. \qed
\medskip

\noindent{\em Proof of Theorem~\ref{general}:}
As the statement of the theorem is asymptotic in nature we may assume that $|T|/n$ is sufficiently large. Note that $|T|/n<n$, so this also means that $n$ is sufficiently large. We introduce the notation $f=|X|/|T|$.

We first reduce the maximum number of intersection points on a curve in $\A$ to $d=\lfloor |X|/n\rfloor$.
To this end, we break each arc $a\in \A$ into sub-arcs $a_1,\ldots, a_h$ all of which, with the possible exception of the last one, contain exactly $d$ intersection points. This splitting yields a family $\A'$ of $n'=\Theta(n)$ Jordan arcs without modifying the set $X$ of intersection points or the set $T$ of touching points.

If $|T|^2<n'd^3$ we apply Lemma~\ref{denselemma} to $\A'$. The lemma claims that $f^{72}h^{144}=\Omega(\log\log n')$ for $h=n^{\prime2}/|T|$. We further have $|T|^2<n'd^3=O(|X|^3/n^2)$ yielding $h<O(f^3)$ and thus $f^{72}h^{144}=O(f^{504})$. The statement of the theorem follows.

If $|T|^2\ge n'd^3$, then we apply Lemma~\ref{densepart} to the collection $\A'$ and the set of touching points $T$. Let $\A'_0\subseteq\A'$ be the collection whose existence is claimed by this lemma. The size $|\A'_0|$ of this collection is $n'_0=\Theta(n'^2d^3/|T|^2)=\Theta(f^2|X|/n)$. The family determines $x'_0\le n'_0d=O(f^2|X|^2/n^2)$ intersections among which $t'_0=\Omega(n'd^3/|T|)=\Omega(f|X|^2/n^2|)$ are touchings. We apply Lemma~\ref{denselemma} to the family $\A'_0$. With $f'_0=x'_0/t'_0=O(f)$ and $h'_0=n'^2_0/t'_0=O(f^3)$, the lemma states $f'^{72}_0h'^{144}_0=\Omega(\log\log n'_0)$. Here $f'^{72}_0h'^{144}_0=O(f^{504})$ and $n'_0=\theta(f^2|X|/n)=\Omega(|X|/n)=\Omega(|T|/n)$. The statement of the theorem follows again. \qed


\begin{thebibliography}{}


\bibitem[ANPPSS04]{PseudoCircles} P. K. Agarwal, E. Nevo, J. Pach, R. Pinchasi, M. Sharir, and S. Smorodinsky, Lenses in arrangements of pseudocircles and their applications, {\it J. ACM} {\bf 51} (2004), 139--186.



\bibitem[AgS05]{AgS05} P. K.  Agarwal and M. Sharir, Pseudo-line arrangements: duality, algorithms, and applications, {\em  SIAM J. Comput.} {\bf 34} (2005), no. 3, 526--552.

\bibitem[ACNS82]{ACNS} M. Ajtai, V. Chv\'atal, M. Newborn, E. Szemer\'edi, Crossing-free subgraphs. In: {\em Theory and Practice of Combinatorics}, North-Holland Mathematics Studies {\bf60}, North-Holland, Amsterdam, 1982, pp. 9–-12.

\bibitem[An70]{An70} E. M. Andreev, Convex polyhedra of finite volume in Loba\v{c}evski\u{i} space, {\it Mat. Sb. (N.S.)}, {\bf 83} (1970), no. 125, 256 –- 260.

\bibitem[ArS02]{ArS02}  B. Aronov and M. Sharir, Cutting circles into pseudo-segments and improved bounds for incidences, {\em Discrete Comput. Geom.} {\bf 28} (2002), no. 4, 475--490.


\bibitem[Ch1]{Chan1} T. M. Chan, On levels in arrangements of curves, Discrete Comput. Geom. {\bf 29} (2003), 375--393. (Also in {\it Proc. 41th IEEE Sympos. Found. Comput. Sci. (FOCS)}, 2000, pp. 219--227.)

\bibitem[Ch2]{Chan2} T. M. Chan, On levels in arrangements of curves, II: a simple inequality and its consequence, {\it Discrete Comput. Geom.} {\bf 34} (2005), 11--24. (Also in {\it Proc. 44th IEEE Sympos.  Found. Comput. Sci. (FOCS)}, 2003, pp. 544--550.)

\bibitem[Ch3]{Chan3} T. M. Chan, On levels in arrangements of curves, III: further improvements, {\it Proc. 24th ACM Symposium on Computational Geometry (SoCG)}, 2008, pp. 85--93.

\bibitem[Cha00]{Cha00} B.  Chazelle, {\em The Discrepancy Method. Randomness and Complexity}, Cambridge University Press, Cambridge, 2000.

\bibitem[CEGSW90]{CEGSW90} K. L. Clarkson, H. Edelsbrunner, L. J. Guibas, M. Sharir, and E. Welzl, Combinatorial complexity bounds for arrangements of curves and spheres, {\em Discrete Comput. Geom.} {\bf 5} (1990), no. 2, 99--160.

\bibitem[CS89]{CS89} K. L. Clarkson and P. W. Shor, Applications of random sampling in computational geometry. II, {\em Discrete Comput. Geom.} {\bf 4} (1989), no. 5, 387--421.

\bibitem[De98]{De98} T. K. Dey, Improved bounds for planar $k$-sets and related problems, {\em Discrete Comput. Geom.} {\bf 19} (1998), no. 3, 373--382.

\bibitem[Dv10]{Dv10} Z. Dvir, Incidence theorems and their applications, {\em Found. Trends Theor. Comput. Sci.} {\bf 6} (2010), no. 4, 257--393 (2012).

\bibitem[Ed87]{Ed87} H. Edelsbrunner, {\em Algorithms in Combinatorial Geometry. EATCS Monographs on Theoretical Computer Science, 10}, Springer-Verlag, Berlin, 1987.

\bibitem[ESZ16]{ESZ16} J. S. Ellenberg, J. Solymosi, and J. Zahl, New bounds on curve tangencies and orthogonalities, {\it Disc. Anal.} {\bf 22} (2016), 1 -- 22.

\bibitem[EH89]{EH89} P. Erd\H os and A. Hajnal, Ramsey-type theorems, {\em Discrete Appl. Math.} {\bf 25} (1989), no. 1-2, 37--52.

\bibitem[Er46]{Er46} P. Erd\H os, On sets of distances of n points, {\em Amer. Math. Monthly} {\bf 53} (1946), 248--250.

\bibitem[FFPP10]{FFPP10}  J. Fox, F. Frati, J. Pach, and R. Pinchasi: Crossings between curves with many tangencies, in: {\em WALCOM: Algorithms and Computation, Lecture Notes in Comput. Sci.} {\bf 5942}, Springer-Verlag, Berlin, 2010, 1--8. Also in: {\em An Irregular Mind, Bolyai Soc. Math. Stud.} {\bf 21}, J\'anos Bolyai Math. Soc., Budapest, 2010, 251--260.


\bibitem[FP08]{FP08} J. Fox and J. Pach, Separator theorems and Tur\'{a}n-type results for planar intersection graphs, {\it Adv.
in Math.} 219 (2008), 1070--1080.

\bibitem[FP10]{FP10} J. Fox and J. Pach, A separator theorem for string
graphs and its applications, {\it Combin. Probab. Comput.}
19 (2010), 371--390.

\bibitem[FP12]{FP12} J. Fox and J. Pach, String graphs and incomparability
graphs, {\it Adv. in Math.} 230 (2012), 1381--1401.

\bibitem[FPT11]{FPT11} J. Fox, J. Pach, and C. D. T\'{o}th, Intersection patterns
of curves, {\it J. London Math. Soc.} 83 (2011), 389-406.



\bibitem[G]{G3} L. Guth, Polynomial partitioning for a set of varieties, manuscript.


\bibitem[GK10]{GK10} L. Guth and N. H. Katz, Algebraic methods in discrete analogs of the Kakeya problem, {\em Adv. Math.} {\bf 225} (2010), no. 5, 2828--2839.

\bibitem[G13]{G1} L. Guth, Unexpected applications of polynomials in combinatorics, in: {\em The Mathematics of Paul Erd\H os. I,} Springer, New York, 2013, 493--522.

\bibitem[G15]{G2} L. Guth, Distinct distance estimates and low degree polynomial partitioning, {\em Discrete Comput. Geom.} {\bf 53} (2015), no. 2, 428--444.

\bibitem[GK15]{GK15} L. Guth and N. H. Katz, On the Erd\H os distinct distances problem in the plane, {\em Ann. of Math. (2)} {\bf 181} (2015), no. 1, 155--190.

\bibitem[HW87]{HW87} D. Haussler and E. Welzl, $\varepsilon$-nets and simplex range queries, {\em Discrete Comput. Geom.} {\bf 2} (1987), no. 2, 127--151.

\bibitem[KLPS86]{KLPS86} K. Kedem, R. Livn\'{e}, J. Pach, and M. Sharir, On the union of Jordan regions and collision-free translational motion amidst polygonal obstacles, {\it Discrete Comput. Geom.} {\bf 1 (1)} (1986), 59 -- 71.


\bibitem[Koe36]{Koe36} P. Koebe, Kontaktprobleme der Konformen Abbildung, {\it Ber. S\"{a}chs. Akad. Wiss. Leipzig, Math.-Phys. Kl.} {\bf 88} (1936), 141 -– 164.

\bibitem[KST54]{KST54} T. K\H{o}v\'{a}ri, V. S\'{o}s, and P. Tur\'{a}n, On a problem of K. Zarankiewicz, {\it Colloquium Math.}
{\bf 3} (1954), 50 -– 57.

\bibitem[Le83]{L} T. Leighton, Complexity Issues in VLSI, {\em Foundations of Computing Series}, MIT Press, Cambridge, 1983.

\bibitem[MaT06]{MaT06}  A. Marcus and G. Tardos, Intersection reverse sequences and geometric applications, {\em J. Combin. Theory Ser. A} {\bf 113} (2006), no. 4, 675--691.


\bibitem[Ma14]{Ma14} J. Matou\v sek, Near-optimal separators in string graphs, {\em Combin. Probab. Comput.} {\bf 23} (2014), no. 1, 135--139.

\bibitem[Mu02]{Mu02}  D.  Mubayi: Intersecting curves in the plane, {\em Graphs Combin.} {\bf 18} (2002), no. 3, 583--589.


\bibitem[PRT15]{PRT15} J. Pach, N. Rubin, and G. Tardos, On the Richter-Thomassen conjecture about pairwise intersecting closed curves, in: {\em Proc. 26th Annual ACM-SIAM Symposium on Discrete Algorithms} (SODA 2015, San Diego), SIAM, 2015, 1506--1515. Also {\it Combin. Prob. Comput.} {\bf 25} (2016), no. 6, 941 -- 958.

\bibitem[PRT16]{PRT16}J. Pach, N. Rubin, and G. Tardos, Beyond the Richter-Thomassen Conjecture, in: {\em Proc. 27th Annual ACM-SIAM Symposium on Discrete Algorithms} (SODA 2016, Arlington), SIAM, 2016, 957--968.

\bibitem[PaS98]{PaS98}  J. Pach and M. Sharir, On the number of incidences between points and curves, {\em Combin. Probab. Comput.} {\bf 7} (1998), no. 1, 121--127.

\bibitem[PaS09]{PaS09} J. Pach and M. Sharir, {\em Combinatorial Geometry and its Algorithmic Applications. The Alcal\'a lectures. Mathematical Surveys and Monographs, 152}, American Mathematical Society, Providence, RI, 2009.

\bibitem[PST12]{PST12} J. Pach, A. Suk, and M. Treml, Tangencies between families of disjoint regions in the plane, {\em Comput. Geom.} {\bf 45} (2012), no. 3, 131--138.

\bibitem[RiT95]{RiT95}  R. B. Richter and C. Thomassen, Intersections of curve systems and the crossing number of $C_5\times C_5$, {\em Discrete Comput. Geom.} {\bf 13} (1995), no. 2, 149--159.

\bibitem[Sa99]{Sa99}  G. Salazar: On the intersections of systems of curves, {\em J. Combin. Theory Ser. B} {\bf 75} (1999), no. 1, 56--60.

\bibitem[ShA95]{ShA95}  M. Sharir and P. K. Agarwal: {\em Davenport-Schinzel Sequences and Their Geometric Applications}, Cambridge University Press, Cambridge, 1995.

\bibitem[SSZ15]{SSZ15} M. Sharir, A. Sheffer, and J. Zahl, Improved bounds for incidences between points and circles, {\em Combin. Probab. Comput.} {\bf 24} (2015), no. 3, 490--520.



\bibitem[SzT83a]{SzT83a} E. Szemer\'edi and W. T. Trotter, Jr., Extremal problems in discrete geometry, {\em Combinatorica} {\bf 3} (1983), no. 3-4, 381--392.

\bibitem[SzT83b]{SzT83b} E. Szemer\'edi and W. T. Trotter, Jr., A combinatorial distinction between the Euclidean and projective planes, {\em European J. Combin.} {\bf 4} (1983), no. 4, 385--394.

\bibitem[TT98]{TT98} H. Tamaki and T. Tokuyama, How to cut pseudoparabolas into segments, {\em Discrete Comput. Geom.} {\bf 19} (1998), no. 2, 265--290.

\bibitem[Thu97]{Thu97} W. Thurston, Three-dimensional Geometry and Topology, Vol. 1, Silvio Levy Ed., Princeton Mathematical Series, 35, Princeton University Press, Princeton, NJ, 1997.


\bibitem[Tut70]{Tut70} W. T. Tutte, Toward a theory of crossing numbers, J. Combinatorial Theory 8 (1970), 45 -– 53.



\end{thebibliography}
\end{document}